\documentclass[twoside,10pt]{article}
\usepackage{etex}
\usepackage{lipsum} 

\usepackage[sc]{mathpazo} 
\usepackage[T1]{fontenc} 
\usepackage{microtype} 

\usepackage[hmarginratio=1:1,top=32mm,columnsep=20pt]{geometry} 
\usepackage{multicol} 
\usepackage[hang, small,labelfont=bf,up,textfont=it,up]{caption} 
\usepackage{booktabs} 
\usepackage{paralist} 
\usepackage{amssymb}
\usepackage{abstract} 

\usepackage{titlesec} 
\titleformat{\section}[block]{\large\scshape}{\thesection.}{1em}{} 
\titleformat{\subsection}[block]{\large}{\thesubsection.}{1em}{} 

\usepackage{fancyhdr} 
\usepackage{amsthm,floatrow}
\usepackage{amsmath,amstext,graphicx,epstopdf}
\usepackage[position=t,singlelinecheck=off]{subfig}
\usepackage{multirow,tikz}
\usetikzlibrary{calc}
   \usepackage[yyyymmdd,hhmmss]{datetime}
\usepackage{tipa}
\usepackage{caption}
\usepackage{algorithmic}
\usepackage[ruled,vlined]{algorithm2e}
\usepackage{bm}
\usepackage{adjustbox}
   
\makeatletter
\newcounter{megaalgorithm}

\makeatother

\makeatletter
\newcounter{algorithmone}

\makeatother

\makeatletter
\newcounter{algorithmtwo}

\makeatother

\newtheorem{dummy}{***}

\newtheorem{defn}[dummy]{Definition}

\usepackage{xparse}

\NewDocumentCommand{\ceil}{s O{} m}{%
  \IfBooleanTF{#1} 
    {\left\lceil#3\right\rceil} 
    {#2\lceil#3#2\rceil} 
}

\usepackage{parskip}
\usepackage{lipsum}
\makeatletter
\def\thm@space@setup{%
  \thm@preskip=12pt plus 0pt minus 8pt}
\makeatother

\pgfmathsetmacro{\minsize}{0.2cm}

\newlength\myindent
\newcommand\bindent{%
  \begingroup
  \setlength{\itemindent}{\myindent}
  \addtolength{\algorithmicindent}{\myindent}
}
\newcommand\eindent{\endgroup}

\newlength\myindenttwo
\newcommand\bindenttwo{%
  \begingroup
  \setlength{\itemindent}{\myindenttwo}
  \addtolength{\algorithmicindent}{\myindenttwo}
}
\newcommand\eindenttwo{\endgroup}

\newlength\myindentthree
\newcommand\bindentthree{%
  \begingroup
  \setlength{\itemindent}{\myindentthree}
  \addtolength{\algorithmicindent}{\myindentthree}
}
\newcommand\eindentthree{\endgroup}

\usepackage{array}
\newcolumntype{M}[1]{>{\centering\arraybackslash}m{#1}}
\newcolumntype{C}[1]{>{\centering\arraybackslash$}p{#1}<{$}}

\newcommand{\diagdots}[3][-25]{%
  \rotatebox{#1}{\makebox[0pt]{\makebox[#2]{\xleaders\hbox{$\cdot$\hskip#3}\hfill\kern0pt}}}%
}

\newcommand{\tstar}[7]{
\pgfmathsetmacro{\starangle}{360/#3}
\begin{scope}[shift={(#6+0.5,#7+0.5)}]
\draw[#5] (#4:#1)
\foreach \x in {1,...,#3}
{ -- (#4+\x*\starangle-\starangle/2:#2) -- (#4+\x*\starangle:#1)
}
-- cycle ;
\end{scope}
}

\newcommand{\ngram}[6]{
\pgfmathsetmacro{\starangle}{360/#2}
\pgfmathsetmacro{\innerradius}{#1*sin(90-\starangle)/sin(90+\starangle/2)}
\tstar{\innerradius}{#1}{#2}{#3}{#4}{#5}{#6}
}

\title{\vspace{-15mm}\fontsize{16pt}{10pt}\selectfont
\textbf{A Multigrid Algorithm with Non-Standard Smoother for Two Selective Models in Variational Segmentation}}
\author{
\large
\textsc{Michael Roberts$^\dagger$, Ke Chen$^\dagger$\thanks{Email k.chen@liverpool.ac.uk,\ Web: www.liv.ac.uk/cmit \ (corresponding author). Work supported by UK EPSRC   grant EP/K036939/1.}\,  and
Klaus Irion$^\ddagger$
}
\\[2mm] 
\normalsize $^\dagger$Centre for Mathematical Imaging Techniques and Department of Mathematical\\
\normalsize   Sciences, The University of Liverpool,
\normalsize   United Kingdom\\
and \normalsize $^\ddagger$Department of Radiology, Liverpool Heart and Chest Hospital, \\
\normalsize  Liverpool, United Kingdom
}
\vspace{-5mm}
\date{}

\begin{document}

\maketitle


\begin{abstract}
\small
Automatic segmentation of an image to identify all meaningful parts is one of the most challenging as well as useful tasks in
a number of application areas. This is widely studied. Selective segmentation, less studied, aims to use limited user specified information to
extract one or more interesting objects (instead of all objects).
Constructing a fast solver remains a challenge for both classes of model.
However our primary concern is on selective segmentation.

In this work, we develop an effective multigrid algorithm, based on a new non-standard smoother to deal with non-smooth coefficients, to solve the underlying partial differential equations (PDEs) of a class of variational segmentation models  in the level set formulation.
For such models, non-smoothness (or jumps) is typical as segmentation is only possible if edges (jumps) are present.
In comparison with previous multigrid methods   which were shown to produce an acceptable {\it mean} smoothing rate for related models, the new algorithm can ensure a small and {\it global} smoothing rate that is a sufficient condition for convergence.
Our rate analysis is by Local Fourier Analysis and, with it, we design the corresponding iterative solver, improving on an ineffective line smoother. Numerical tests show that the new algorithm outperforms multigrid methods based on competing smoothers.


\end{abstract}


\section{Introduction}

Segmentation of an image into its individual objects is one incredibly important
application of image processing techniques. Not only are accurate segmentation results
required, but also it is required that the segmentation method is fast.
Many imaging applications demand increasingly higher resolution e.g.
an image of size 25000 $\times$ 25000 (or practically $10^8$  unknowns)
can be common in oncology imaging.
Here we address the problem of slow solutions by developing a
 fast multigrid method for PDEs arising from segmentation models.


Segmentation can take two forms; firstly global segmentation is the isolation of all objects
 in an image from the background and secondly, selective segmentation is the
 isolation of a subset of the objects in an image from the background. Selective
 segmentation is very useful in, for example, medical imaging for the segmentation
 of single organs.

Approaches to image segmentation broadly fall into two classes; region-based and edge-based. Some region-based approaches are region growing \cite{RC1}, watershed algorithms \cite{BC31}, Mumford-Shah \cite{SC22} and Chan-Vese \cite{SC10}. The final two of these are PDE-based variational approaches to the problem of segmentation.
There are also models which mix the two classes to use the benefits of the region-based and edge-based approaches and will incorporate features of each.
Edge-based methods aim to encourage an evolving contour towards the edges in an image
and normally require
 an edge detector function \cite{SC5}. 
The first edge-based variational approach was devised by Kass et al.   \cite{SC17}
with the famous snakes model, this was further developed by Casselles et al.  \cite{SC5}
who introduced the Geodesic Active Contour (GAC) model.
Region-based global segmentation models include the well known works of
 Mumford-Shah \cite{SC22} and Chan-Vese \cite{SC10}.
 Importantly they are non-convex and hence a minimiser of these models may only be a local, not the global, minimum. Further work by Chan et al. \cite{SC8} gave rise to a method to find
the global minimiser for the Chan-Vese model under certain conditions.

Selective segmentation of objects in an image, given a set of points near the object or objects to be segmented,
 builds in such user input to a model
using a set $
\mathcal{S}=\{ (x_{i},y_{i})\in\Omega, 1\le i\le k\}
$
where $\Omega\subset\mathbb{R}^{2}$ is the image domain \cite{RC16, SC1, SC2}.
 Nguyen et al. \cite{SC25}  considered marker sets $\mathcal{S}$ and $\mathcal{A}$ which consist of points inside and outside, respectively, the object or objects to be segmented.
Gout et al. \cite{RC16} combined the GAC approach with the geometrical constraint that
 the contour pass through the points of $\mathcal{S}$. This was enforced with a
 distance function which is zero at $\mathcal{S}$ and non-zero elsewhere.
Badshah and Chen \cite{SC1} then   combined the Gout et al. model with
  \cite{SC10}  to incorporate a constraint on the intensity in the
selected region, thereby encouraging the contour to segment homogenous regions.
Rada and Chen \cite{SC27} introduced a selective segmentation method based on
 two-level sets
which was shown to be more robust than the Badshah-Chen model.
We also refer to \cite{SC2,MR1} for selective segmentation models which include
different fitting constraints, using coefficient of variation and
the centroid of $\mathcal{S}$ respectively.

None of these models have a restriction on the size of the object or objects to be detected and depending on the initialisation these methods have the potential to detect more or
 fewer objects than the user desired.
To address this and to improve on \cite{SC27}, Rada and Chen \cite{SC28} introduced a
model (we refer to it as the Rada-Chen model from now on) combining the Badshah-Chen \cite{SC1} model with a constraint on the area of
the objects to be segmented.
The reference area used to constrain the area within the
contour is that of the polygon formed by the markers in $\mathcal{S}$.
Spencer and Chen \cite{SC} recently introduced a model with the distance fitting
penalty as a standalone term in the energy functional, unbounding it from the edge
detector term of the Gout et al. model.
All of the above selective segmentation models discussed are non-convex and hence
the final result depends on the initialisation. Spencer and Chen \cite{SC}, in the same paper,
reformulated the model they introduced to a convex form using a penalty term as in \cite{SC8}. We have considered the convex Spencer-Chen model but
found that the numerical implementation is unfortunately
 sensitive to the main parameters and is unstable if they aren't chosen correctly within a small range; hence we focus on the non-convex
model they introduce for which reliable results have been found (we refer to this as the
Spencer-Chen model from now on).
A convex version of the Rada-Chen model cannot be formulated \cite{SC}.
In this paper we only consider 2D images, however for completion we remark that 3D segmentation models do exist \cite{SC13,SC31}.

Solving the PDE models, in the context of large scale images,
quickly remains a challenge.
The variational approach to image segmentation involves the minimisation of an energy functional such as that in \cite{SC28}. This will typically involve solving a system of equations from a discretised PDE using an iterative method. In particular, discretisations of models such as \cite{SC1, SC2, SC10, SC28, SC} are non-linear and so require non-linear iterative methods to solve.
The number of equations in the system is equal to the number of pixels in the image, which can be very large, and for each equation in the system the number of steps of an iterative method required can also be very large (to reach convergence).
Due to improvements in technology and imaging, we now can produce larger and larger images, however this has the direct consequence that analysis of such images has become much more computationally intensive. We remark that if we directly discretise the variational models first (without using PDEs),
  Chan-Vese type models can be reformulated into minimisation based on graph cuts and then fast algorithms have been proposed \cite{Tai09, Lu16}.

The multigrid approach for solving PDEs in imaging has been tried before
and previous work by Badshah and Chen \cite{BC1,BC2} introduced a 2D Chan-Vese multigrid algorithm for two-phase and multi-phase images, additionally Zhang et al. \cite{SC31} implemented a multigrid algorithm for the 3D Chan-Vese model.
The fundamental idea behind multigrid is that if we perform most of the computations on a reduced resolution image then the computational expense is lower.
We then transfer our solution from the low resolution grid to the high resolution grid through interpolation and smooth out any errors which have been introduced by the interpolation using a few steps of a smoothing algorithm, e.g. Gauss-Seidel.
The multigrid method is an optimal solver when it converges \cite{Lin:16,BC30}. This requires that the smoothing scheme, which corrects the errors when transferring between the higher and lower resolution images
and vice-versa, is effective, i.e. reduces the error magnitude
of high-frequency components quickly.

In the large literature of multigrid methods, the convergence problem associated with non-smooth or jumping coefficients was often highlighted \cite{CB1,CB} and developing working algorithms which converge is a key problem. Much attention was given to
designing better coarsening strategies and improved interpolation operators
\cite{Wan10, Zhu14} while keeping the simple smoothers; such as the damped Jacobi,
Gauss-Seidel or line smoothers.
In practice, one can quickly exhaust the list of standard smoothers and yet cannot find
a suitable one unless compromising in optimality by increasing the number of iterations.
In contrast, our approach here is to seek a non-standard and
more effective smoother with an acceptable smoothing rate. Our work is motivated by
 Napov and Notay \cite{Napov2011} who established the explicit relationship
 of a smoothing rate to the underlying multigrid convergence rate for linear models; in particular
 the former
 also serves as the lower bound for the latter.

The contributions of this paper can be summarised as follows: (1) We review six smoothers for the Rada-Chen and Spencer-Chen selective segmentation models and perform Local Fourier Analysis (LFA) to assess their performance and quantitatively determine their effectiveness (or lack of).
(2) We propose an effective non-linear multigrid method to solve the Rada-Chen model \cite{SC28} and the Spencer-Chen model \cite{SC}, based on a new smoothers that 
add non-standard smoothing steps locally at coefficient jumps.
We recommend in particular one of our new hybrid smoothers which achieves a better smoothing rate than the other smoothers studied and thus gives rise to a multigrid framework which converges to the energy minimiser faster than when standard smoothers are used.

The remainder of this paper is structured as follows; in \S 2 we review some
global and selective segmentation models
building to the Rada-Chen and Spencer-Chen models.
In \S 3 we describe the Full Approximation
Scheme multigrid framework, 
  give details of six smoothers that we consider 
and compare the smoothing rates.
We find that none of these standard smoothers can produce a small enough
smoothing rate to yield an effective  multigrid method and so in \S 4 we then introduce
two new hybrid smoothers based on
new iterative schemes to
improve the smoothing rates at those pixels where the six smoothers perform badly.
 In \S 5 we test our algorithms with
 some numerical results, recommend the best algorithm using one of our proposed smoothers
and analyse the complexity of the recommended multigrid algorithm.
Finally in \S 6 we provide some concluding remarks.


\section{Review of segmentation models}
Our methods will apply to both global segmentation models
and selective segmentation models.
It is necessary to briefly describe both types.
Denote  a given image in domain $\Omega\subset\mathbb{R}^{2}$ by $z(x,y)$.

\subsection{Global segmentation models}

The model of Mumford and Shah \cite{SC22} is one of the most famous and
important variational models in image segmentation.
We will review its two-dimensional piecewise constant variant,
commonly  known as the Chan-Vese (CV) model \cite{SC10},
which takes the form
\begin{equation}\label{eqn:cv}
\min_{\Gamma,c_{1},c_{2}}F_{CV}(\Gamma,c_{1},c_{2})
=\mu\cdot length (\Gamma)+\lambda_{1}\int_{\Omega_{1}}
|z(x,y)-c_{1}|^{2}d\Omega +\lambda_{2}\int_{\Omega_{2}}|z(x,y)-c_{2}|^{2}d\Omega
\end{equation}
where
the foreground $\Omega_{1}$ is the subdomain to be segmented,
  the background is $\Omega_{2}=\Omega\backslash\Omega_{1}$
  and $\mu,\lambda_{1},\lambda_{2}$ are fixed non-negative parameters.
The values $c_{1}$ and $c_{2}$ are the average intensities of $z(x,y)$
inside   $\Omega_{1}$ and $\Omega_{2}$ respectively.
Using the ideas of Osher and Sethian \cite{OS1}, a level set function
\[
\phi(x,y)=
\begin{cases}
    >0,& (x,y)\in\Omega_{1},\\
    0,& (x,y)\in\Gamma,\\
    <0, & otherwise,\\
\end{cases}
\]
is used by \cite{SC10} to   track the object boundary $\Gamma$, where we now define it as the zero level set of $\phi$, i.e.  $\Gamma = \{ (x,y)\in\Omega\, |\, \phi(x,y)=0\}$. We reformulate  (\ref{eqn:cv}) as
\begin{equation}\label{eqn:cvls}
\begin{split}
\min_{\phi,c_{1},c_{2}}F_{CV}(\phi,c_{1},c_{2})=&\mu\int_{\Omega}|\nabla H_{\varepsilon}(\phi)|d\Omega+\lambda_{1}\int_{\Omega}(z(x,y)-c_{1})^{2}H_{\varepsilon}(\phi)d\Omega\\
&
\ \hspace*{2.7cm}+\lambda_{2}\int_{\Omega}(z(x,y)-c_{2})^{2}(1-H_{\varepsilon}(\phi))d\Omega,\\
\end{split}
\end{equation}
with $H_{\varepsilon}(\phi)$ a smoothed Heaviside function such as \cite{SC10}
\[
H_{\varepsilon}(\phi) = \frac{1}{2}+\frac{1}{\pi}\arctan\bigg(\frac{\phi}{\varepsilon}\bigg)
\]
where we use $\varepsilon =1$ in our experiments. We solve this minimisation problem in two stages, first
%
with $\phi$ fixed we minimise with respect to $c_{1}$ and $c_{2}$, yielding
\begin{equation}\label{eqn:c1c2}
c_{1}=\frac{\int_{\Omega}H_{\varepsilon}(\phi)\cdot z(x,y)\, d\Omega}{\int_{\Omega}H_{\varepsilon}(\phi) \,d\Omega},
\hspace{0.5in}
c_{2}=\frac{\int_{\Omega}(1-H_{\varepsilon}(\phi))\cdot z(x,y)\, d\Omega}{\int_{\Omega}(1-H_{\varepsilon}(\phi)) \,d\Omega},
\end{equation}
and secondly, with $c_{1}$ and $c_{2}$ fixed we minimise (\ref{eqn:cvls}) with respect to $\phi$. This requires the determination of the associated Euler-Lagrange form \cite{SC10} and then solving the resulting PDE.
A drawback of the Chan-Vese functional (\ref{eqn:cvls}) is that it is non-convex.
Therefore a minimiser of this functional may only be a local minimum and not the global minimum.
Hence the final segmentation result is dependent on the initial contour. Chan et al. \cite{SC8} reformulated
(\ref{eqn:cvls}) to obtain an equivalent convex model and hence we can always obtain the global minimum for this model.

\subsection{Selective segmentation models}

\label{sec:selsegintro}

Selective segmentation models make use of user input, being a marker set of points near the object or objects to be segmented. Let
\(
\mathcal{S}=\{ (x_{i},y_{i})\in\Omega, 1\le i\le k\}
\)
be such a marker set. The contour is encouraged to pass through or near the points of $\mathcal{S}$ by a distance function such as \cite{SC13}
\[
d(x,y)=\prod^{k}_{i=1}\Big( 1- e^{-\frac{(x_{i}-x)^{2}}{2\sigma^{2}}}e^{-\frac{(y_{i}-y)^{2}}
{2\sigma^{2}}}\Big),\ \ \forall (x,y)\in\Omega, (x_{i},y_{i})\in\mathcal{S},
\]
where $\sigma$ is a fixed non-negative tuning parameter. See, for example, \cite{RC16, SC} for other distance functions.
The distance function is zero at the points of $\mathcal{S}$ and non-zero elsewhere,
taking a maximum value of one. Gout et al. \cite{SC13} were the first to introduce
a model incorporating a distance function into the Geodesic Active Contour model
of Caselles et al. \cite{SC5}, however
this model struggles when boundaries between objects and their background
are fuzzy or blurred. 
To address this, Badshah and Chen \cite{SC1} introduced a new model which includes
 the intensity fitting terms from the CV model (\ref{eqn:cv}).
%
However this model has poor robustness \cite{SC27} if iterating for too many steps
the final segmentation can include  more or fewer objects than intended.
To improve on this, Rada and Chen \cite{SC28} introduced a model which incorporates
an area fitting term into the Badshah-Chen (BC) model and is far more robust.

{\bf The Rada-Chen model} \cite{SC28}.
This is the first model we focus on in this paper, defined by
\begin{equation}\label{eqn:rcfunc}
\begin{split}
F_{RC}(\phi,c_{1},c_{2})=&\mu\int_{\Omega}d(x,y)g(|\nabla z(x,y)|^{2})|\nabla
H_{\varepsilon}(\phi)|dxdy \\
& + \lambda_{1}\int_{\Omega}(z(x,y)-c_{1})^{2}H_{\varepsilon}(\phi)dxdy +
\lambda_{2}\int_{\Omega}(z(x,y)-c_{2})^{2}(1-H_{\varepsilon}(\phi))dxdy\\
& + \nu\bigg[\left(\int_{\Omega}H_{\varepsilon}(\phi)dxdy-A_{1}\right)^{2}+
\left(\int_{\Omega}(1-H_{\varepsilon}(\phi))dxdy-A_{2}\right)^{2}\bigg],
\end{split}
\end{equation}
where $\mu,\lambda_{1},\lambda_{2},\nu$ are fixed non-negative parameters. The edge detector function $g(|\nabla z(x,y)|^{2})$ is given by $g(s) = 1/(1+\beta s)$ for tuning parameter $\beta$ which takes value $0$ at edges and is $1$ away from them. $A_{1}$ is the area of the polygon formed from the points of $\mathcal{S}$ and $A_{2}=|\Omega|-A_{1}$.
The final term of this functional therefore puts a penalty on the area inside a contour being very different to $A_{1}$.
The first variation of (\ref{eqn:rcfunc}) with respect to $\phi$ gives the Euler-Lagrange form   \cite{SC28}
\begin{equation}\label{eqn:rcel}
\begin{gathered}
\delta_{\varepsilon}(\phi)\bigg\{\mu\nabla\cdot\left(\frac{d(x,y)\cdot g(|\nabla z(x,y)|^{2})
\nabla\phi}{|\nabla\phi|}\right)-\Big[\lambda_{1}(z(x,y)-c_{1})^{2}-\lambda_{2}
(z(x,y)-c_{2})^{2}\Big]\bigg.\\
\bigg.-\nu\bigg[(\int_{\Omega}H_{\varepsilon}(\phi)dxdy-A_{1})-(\int_{\Omega}
(1-H_{\varepsilon}(\phi))-A_{2})\bigg]\bigg\}=0,\\
\end{gathered}
\end{equation}
in $\Omega$ with the condition that $\frac{\partial\phi}{\partial\bm{n}}=0$ on $\partial\Omega$, $\bm{n}$ the outward normal vector and $\delta_{\varepsilon}(\phi) = \frac{d H_{\varepsilon}(\phi)}{d\phi}$.

{\bf Discretisation of the Rada-Chen model}.
We denote by $\phi_{i,j}=\phi(x_{i},y_{j})$ the approximation of
$\phi$ at $(i,j)$ for $1\le i\le n$ and $1\le j\le m$.
We let $h_{x}$ and $h_{y}$ be the grid spacings in the $x$ and $y$ directions respectively. Using finite differences, and noting $A_{2}=1-A_{1}$, we obtain the scheme
\begin{equation}\label{eqn:rcabcd}
\begin{gathered}\ \hspace*{-41mm}
A_{i,j}\phi_{i+1,j}+B_{i,j}\phi_{i-1,j}+C_{i,j}\phi_{i,j+1}+D_{i,j}\phi_{i,j-1}-S_{i,j}\phi_{i,j}\\ \qquad
-\delta_{\varepsilon}(\phi_{i,j})\bigg\{\bigg[\lambda_{1}(z_{i,j}-c_{1})^{2}-\lambda_{2}
(z_{i,j}-c_{2})^{2}\bigg]
-2\nu\bigg[h_{x}h_{y} \sum_{k,l}H_{\varepsilon}(\phi_{k,l})-A_{1}\bigg] \bigg\}=0,\\
\end{gathered}
\end{equation}
where\qquad
\(\displaystyle
G_{i,j}=\frac{d_{i,j}\cdot g(|\nabla z_{i,j}|)}{|\nabla\phi_{i,j}|}, \quad
A_{i,j}=\frac{\mu\delta_{\varepsilon}(\phi_{i,j})}{h_{x}^{2}}G_{i+\frac{1}{2},j},
\quad
B_{i,j}=\frac{\mu\delta_{\varepsilon}(\phi_{i,j})}{h_{x}^{2}}G_{i-\frac{1}{2},j},
\)
\begin{equation}\label{eqn:abcdrc}
\begin{gathered}
C_{i,j}=\frac{\mu\delta_{\varepsilon}(\phi_{i,j})}{h_{y}^{2}}G_{i,j+\frac{1}{2}},\quad
D_{i,j}=\frac{\mu\delta_{\varepsilon}(\phi_{i,j})}{h_{y}^{2}}G_{i,j-\frac{1}{2}},
\quad
S_{i,j}=A_{i,j}+B_{i,j}+C_{i,j}+D_{i,j},
\end{gathered}
\end{equation}
%

{\bf The Spencer-Chen model} \cite{SC}.
The second model we focus on in this paper is defined by
\begin{equation}\label{eqn:scfunc}
\begin{split}
F_{SC}(\phi,c_{1},c_{2})=&\mu\int_{\Omega}g(|\nabla z(x,y)|^{2})
|\nabla H_{\varepsilon}(\phi)|dxdy
  + \lambda_{1}\int_{\Omega}(z(x,y)-c_{1})^{2}H_{\varepsilon}(\phi)dxdy\\
  & + \lambda_{2}\int_{\Omega}(z(x,y)-c_{2})^{2}(1-H_{\varepsilon}(\phi))dxdy
  + \theta\int_{\Omega}d(x,y)H_{\varepsilon}(\phi)dxdy,
\end{split}
\end{equation}
where $\mu,\lambda_{1},\lambda_{2}$ and $\theta$ are fixed non-negative parameters.
Note that this model differs from the Rada-Chen model (\ref{eqn:rcfunc}) as
the distance function has been separated from the edge detector term and
is now a standalone penalty term.
This model has Euler-Lagrange form
\begin{equation}\label{eqn:scel}
\begin{split}
\delta_{\varepsilon}(\phi)\left\{\mu\nabla\cdot\left(\frac{g(|\nabla z(x,y)|^{2})\nabla\phi}{|\nabla\phi|}\right)-\left[\lambda_{1}(z(x,y)-c_{1})^{2}-\lambda_{2}(z(x,y)-c_{2})^{2}\right]-\theta d(x,y)\right\}=0,\\
\end{split}
\end{equation}
in $\Omega$ with the condition that $\frac{\partial\phi}{\partial\bm{n}}=0$ on $\partial\Omega$, again with $\bm{n}$ the outward normal vector. We discretise this similarly to the Rada-Chen model previously.

\section{Non-linear multigrid Algorithm 1}

Segmentation using a non-linear multigrid algorithm has been explored by Badshah and Chen \cite{BC1, BC2} for the Chan-Vese model \cite{SC10} and the Vese-Chan model \cite{VC1} which are global segmentation models.
A multigrid method has not yet been applied to selective segmentation and this is the main task of this paper, to apply the multigrid method to the Rada-Chen (\ref{eqn:rcfunc}) and Spencer-Chen (\ref{eqn:scfunc}) selective segmentation models. However as we will see shortly, the task is challenging as standard methods do not work. For brevity we will restrict consideration just to the Rada-Chen model as the derivations for the Spencer-Chen model are similar.

\subsection{The Full Approximation Scheme}\label{sec:fas}

To solve the Rada-Chen model we must solve the non-linear system
(\ref{eqn:rcabcd})
and so we will use the non-linear Full Approximation Scheme \cite{BC9, BC13, BC16, BC30} algorithm due to Brandt \cite{BC3}.
Denote a discretised system by 
\begin{equation}\label{eqn:syseqns}
N^{h}\phi^{h}=f^{h},
\end{equation}
where $h$ indicates that these are the functions on the $n\times m$ cell-centred grid $\Omega^{h}$ and $N^{h}$ is the discretised non-linear operator (which contains the boundary conditions). Similarly define the grids $\Omega^{2h}$ as the $\frac{n}{2}\times\frac{m}{2}$ cell-centred grid resulting from the standard coarsening \cite{BC30} of $\Omega^{h}$, we indicate functions on $\Omega^{2h}$ by $f^{2h},N^{2h}$ and $\phi^{2h}$. Let $\Phi^{h}$ be an approximation to $\phi^{h}$ such that the error $e^{h}=\phi^{h}-\Phi^{h}$ is smooth. Define the residual as $r^{h}=f^{h}-N^{h}\Phi^{h}$. Therefore using (\ref{eqn:syseqns}) we have the residual equation
\[
N^{h}(\Phi^{h}+e^{h})-N^{h}\Phi^{h}=r^{h}.
\]
If the error $e^{h}$ is {\it smooth} then this can be well approximated on $\Omega^{2h}$; the assumption can be a big issue for non-linear problems. With an approximation of $e^{h}$ on $\Omega^{2h}$ we can solve the residual equation on $\Omega^{2h}$, which is significantly less computationally expensive than solving on $\Omega^{h}$, and then transfer this error to $\Omega^{h}$ and use it to correct the approximation $\Phi^{h}$. This method, using the two grids $\Omega^{2h}$ and $\Omega^{h}$, is called a two-grid cycle and it can be nested such that we can consider solving on $\Omega^{4h},\Omega^{8h},\ldots$ and transferring the errors up through the levels to $\Omega^{h}$ and smoothing on each level. This is the multigrid method. We transfer from $\Omega^{h}$ to $\Omega^{2h}$ by restriction and from $\Omega^{2h}$ to $\Omega^{h}$ by interpolation.

{\bf Restriction}.
We use the full-weighting operator $I^{2h}_{h}\Phi^{h}=\Phi^{2h}$ \cite{BC30}
\[
\begin{gathered}
\begin{aligned}
\phi^{2h}_{i,j}=\frac{1}{16}\Big[ \phi^{h}_{2i-1,2j-1}+2\phi^{h}_{2i-1,2j}&+\phi^{h}_{2i-1,2j+1}+2\phi^{h}_{2i,2j-1}+4\phi^{h}_{2i,2j}\\
&+2\phi^{h}_{2i,2j+1}+\phi^{h}_{2i+1,2j-1}+2\phi^{h}_{2i+1,2j}+\phi^{h}_{2i+1,2j+1} \Big],\\%
\end{aligned}
\end{gathered}
\]
and at boundary pixels
$
\phi^{2h}_{i,m}=\frac{1}{2}\left[ \phi^{h}_{2i,m-1}+\phi^{h}_{2i,m}\right]$ and
$\phi^{2h}_{n,j}=\frac{1}{2}\left[ \phi^{h}_{n-1,2j}+\phi^{h}_{n,2j}\right].
$

{\bf Interpolation}.
We use a bilinear interpolation operator $I^{h}_{2h}\Phi^{2h}=\Phi^{h}$ \cite{BC30}
\[
\phi^{h}_{2i,2j}=\phi^{2h}_{i,j},
\hspace{0.5in}
\phi^{h}_{2i+1,2j}=\frac{1}{2}\left[\phi^{2h}_{i,j}+\phi^{2h}_{i+1,j}\right],
\hspace{0.5in}
\phi^{h}_{2i,2j+1}=\frac{1}{2}\left[\phi^{2h}_{i,j}+\phi^{2h}_{i,j+1}\right],
\]
\[
\phi^{h}_{2i+1,2j+1}=\frac{1}{4}\left[ \phi^{2h}_{i,j}+\phi^{2h}_{i+1,j}+\phi^{2h}_{i,j+1}+\phi^{2h}_{i+1,j+1}\right].
\]

We now move to the most important element of the multigrid method -- the smoother. As previously mentioned, we need $e^{h}$ to be smooth to ensure that $\Phi^{h}$ is a good approximation to $\phi^{h}$. In practice, we smooth $e^{h}$ by using an iterative method such as Gauss-Seidel \cite{BC1, BC2} and the success or failure of a multigrid method hinges on the effectiveness of it at smoothing the errors.

\subsection{Smoothers for the Rada-Chen \cite{SC28} model}

Gauss-Seidel and Newton iterative methods have been shown to be effective smoothers
 for PDE problems with smooth coefficients \cite{BC30,Wan10}.
 In this subsection we look at three distinct smoothing iterative techniques;
  lexicographic Gauss-Seidel, line Gauss-Seidel and Newton smoothers.
  For each of these smoothers we consider two different approaches for fixing the
  coefficients in the scheme - globally or locally. Hence overall we consider six smoothers
  for \cite{SC28}; the same smoothers are adaptable for \cite{SC} in a simple way.

{\bf Smoothers 1-2 (GSLEX I - II)}. Lexicographic Gauss-Seidel smoothers are widely used   in
multigrid methods \cite{BC1,BC30}. We update $\phi_{i,j}$ one at a time
and work across and down through the grid of pixels in an image.
{\it Lexicographic Gauss-Seidel smoothers for the Rada-Chen model} \cite{SC28}. We can rearrange (\ref{eqn:rcabcd}) as \begin{equation}\label{eqn:gslexrc}
\phi_{i,j}=\big(A_{i,j}\phi_{i+1,j}+B_{i,j}\phi_{i-1,j}+C_{i,j}\phi_{i,j+1}+D_{i,j}\phi_{i,j-1}
-f_{i,j}\big)\ \big/\ S_{i,j},
\end{equation}
where
$f_{i,j}=\delta_{\varepsilon}(\phi_{i,j})\Big\{\big[\lambda_{1}(z_{i,j}-c_{1})^{2}
-\lambda_{2}(z_{i,j}-c_{2})^{2}\big]
+2\nu\big[
+h_{x}h_{y}\sum_{k,l}H_{\varepsilon}(\phi_{k,l})-A_{1}\big]\Big\},$
to obtain a fixed point scheme for the Rada-Chen model. There are two approaches for implementing this smoother; either update the coefficients globally at the start of each outer iteration or update them locally, immediately after solving for each pixel value. We denote the global smoother by GSLEX-I and the local smoother by GSLEX-II. In the algorithm for both smoothers, we cycle through each pixel $(i,j)$ in turn solving (\ref{eqn:gslexrc}) and updating the value of $\phi(i,j)$, only with GSLEX-II do we update the coefficients immediately and they are used in the update of $\phi(i,j)$ on the next iteration.

{\bf Smoothers 3-4 (GSLINE I - II)}. \label{sec:gsline} Line smoothers are often used for harder problems (e.g. anisotropic coefficients).
Here we perform the Gauss-Seidel updates one column at a time but the approach can be easily reformulated for a row by row update.

{\it Gauss-Seidel line smoothers for the Rada-Chen model} \cite{SC28}.
If we rearrange (\ref{eqn:rcabcd}) to have all the $\phi_{\cdot, j}$ terms on the left hand side we obtain \begin{equation}\label{eqn:gslinerc}
\begin{gathered}
A_{i,j}\phi_{i+1,j}+B_{i,j}\phi_{i-1,j}-S_{i,j}\phi_{i,j} = F_{i,j} = -C_{i,j}\,\phi_{i,j+1}-D_{i,j}\phi_{i,j-1}
+f_{i,j},\\ 
\end{gathered}
\end{equation}
where
we can reformulate (\ref{eqn:gslinerc}) as the following tridiagonal system \begin{equation}\label{eqn:system}
\resizebox{.8\textwidth}{!}{$
\begin{bmatrix}
-S_{1,j} & A_{1,j} & 0 &\hdots & 0 &0 \\ B_{2,j} & -S_{2,j} & A_{2,j} & \ddots& 0 &0\\
0 &B_{3,j} & \ddots & \ddots& \ddots & \vdots \\ \vdots & \ddots & \ddots & \ddots& A_{n-2,j} & 0\\
0 & 0 &\ddots &B_{n-1,j}&-S_{n-1,j}&A_{n-1,j}\\
0 & 0 &\hdots & 0 &B_{n,j} &-S_{n,j}\\
\end{bmatrix}
\cdot
\begin{bmatrix}
\phi_{1,j}\\[2.5mm]
\phi_{2,j}\\
\vdots\\
\vdots\\
\phi_{n-1,j}\\[2.5mm]
\phi_{n,j}\\
\end{bmatrix}=
\begin{bmatrix}
F_{1,j}\\[2.5mm]
F_{2,j}\\
\vdots\\
\vdots\\
F_{n-1,j}\\[2.5mm]
F_{n,j}\\
\end{bmatrix}.
$}
\end{equation}


This system is diagonally dominant (by definition (\ref{eqn:abcdrc})) and if $C_{i,j}+D_{i,j}\neq 0$ then the system is strictly diagonally dominant. We can choose parameters for the edge detector and distance function which ensure this is always true. Therefore this will ensure that the Gauss-Seidel line smoother will converge to a solution \cite{GO1}. As before, we obtain two smoothers; 
the global smoother 
GSLINE-I 
and the local smoother 
GSLINE-II.

{\bf Smoothers 5-6 (NEWT I - II)}. Our last set of smoothers rely on the Newton fixed point iteration schemes.

{\it Newton smoothers for the Rada-Chen model} \cite{SC28}. We can rewrite (\ref{eqn:rcabcd}) in a non-linear form for $\phi_{i,j}$ \[ S_{i,j}\phi_{i,j}^{(k)}-P_{i,j}+Q_{i,j}(\phi_{i,j}^{(k)})=0.
\]
where
$P_{i,j}=A_{i,j}\phi_{i+1,j}+B_{i,j}\phi_{i-1,j}+C_{i,j}\phi_{i,j+1}+D_{i,j}\phi_{i,j-1}-\delta_{\varepsilon}(\phi_{i,j})
\left[\lambda_{1}(z_{i,j}-c_{1})^{2}-\lambda_{2}(z_{i,j}-c_{2})^{2}\right]$ and $Q_{i,j}=2\nu\delta_{\varepsilon}(\phi_{i,j})\left[h_{x}h_{y}
\sum_{k,l}H_{\varepsilon}(\phi_{k,l})-A_{1}\right]$.
The Newton scheme to compute $\phi_{i,j}^{(k+1)}$ is \begin{equation}\label{eqn:newtrc}
\phi_{i,j}^{(k+1)}=\phi_{i,j}^{(k)}-
\big( S_{i,j}\phi_{i,j}^{(k)}-P_{i,j}+ Q_{i,j}(\phi_{i,j}^{(k)})\big) \ \big/ \
\big(S_{i,j}+Q_{i,j}'(\phi_{i,j}^{(k)})\big)
\end{equation}
where $Q_{i,j}'(\phi_{i,j}^{(k)})=2\nu\delta_{\varepsilon}(\phi_{i,j})^{2}h_{x}h_{y}+
2\nu\delta_{\varepsilon}'(\phi_{i,j}))\Big[h_{x}h_{y}\sum_{k,l}
H_{\varepsilon}(\phi_{k,l})-A_{1}\Big]$.
We again have a global smoother, NEWT-I, and a local smoother, NEWT-II.

\subsection{Algorithm 1}
%
In \S\ref{sec:fas} we briefly discussed the FAS across two grids, $\Omega^{h}$ (the fine grid)
and $\Omega^{2h}$ (the coarse grid).
The two-grid cycles 
can be nested so we can perform the majority of the computations on coarser grids than $\Omega^{2h}$,
such as $\Omega^{4h},\Omega^{8h},etc$
and recursive use of $V$-cycles gives rise to multigrid schemes \cite{BC30}.
The general non-linear multigrid Full Approximation Scheme algorithm
is given by Algorithm~\ref{alg:fasscheme}.
\begin{algorithm}[htb!]
\caption{FAS multigrid algorithm, \small{$\phi^{h}\leftarrow FASMG(\phi^{h},N^{h},f^{h},\gamma,\nu_{1},\nu_{2},level,max\_level, Smoother)$}}
\label{alg:nmg}%
\begin{algorithmic}

\STATE{\it\underline{Pre-smoothing:}} Perform $\nu_{1}$ iterations of the smoother:
\bindent
$
\phi^{h}\leftarrow Smoother(\phi^{h},f^{h},\nu_{1}).
$
\eindent

\STATE{\it\underline{Coarse grid correction:}}\,\,\,  Compute the residual:
%
\bindent
$
r^{h}=f^{h}-N^{h}\overline{\phi}^{h}.
$
\eindent

\bindenttwo
\STATE Transfer the residual to $\Omega^{2h}$ by restriction:
$
r^{2h}=I^{2h}_{h}r^{h}.
$
\eindenttwo

\bindenttwo
\STATE Compute:
$
\phi^{2h}=I^{2h}_{h}\phi^{h},
\Phi^{2h}=\phi^{2h},
\overline{f}^{2h}=N^{2h}\phi^{2h}+r^{2h}.
$
\eindenttwo

\bindenttwo
 \IF{$level=max\_level$}
\STATE Compute the exact solution $\phi^{2h}$ of $N^{2h}(\phi^{2h})=N^{2h}(\Phi^{2h})+r^{2h}$
\STATE  on $\Omega^{2h}$ using e.g. time-marching \cite{SC10} or AOS \cite{BC30}.
\ELSE
\STATE Perform $\gamma$ cycles (steps) of
\STATE $\phi^{2h}\leftarrow FASMG(\phi^{2h},N^{2h},f^{2h},\gamma,\nu_{1},\nu_{2},level+1,max\_level,Smoother)$.
\indent \ENDIF
\eindenttwo

\STATE{\it\underline{Interpolation:}}\,\,\,\, Compute:
\bindent
$
e^{2h}=\phi^{2h}-\Phi^{2h}.
$
\eindent

\bindentthree
\STATE Transfer the error to $\Omega^{h}$ by interpolation:
$
e^{h}=I_{2h}^{h}e^{2h}.
$
\eindentthree

\bindentthree
\STATE Correct the fine grid approximation:
$
\phi^{h}=\phi^{h}+e^{h}.
$
\eindentthree

\STATE{\it\underline{Post-smoothing:}} Perform $\nu_{2}$ iterations of the smoother:
\bindent
$
\phi^{h}\leftarrow Smoother(\phi^{h},f^{h},\nu_{2}).
$
\eindent
\end{algorithmic}\label{alg:fasscheme}
\end{algorithm}

\subsection{Local Fourier Analysis of Algorithm 1 for the Rada-Chen Model}

Local Fourier Analysis (LFA) is a useful tool for finding a quantitative measure for the effectiveness of
a smoother \cite{BC3, BC13, BC30}. It is designed to study linear problems with constant coefficients on an infinite grid. However, it is a standard and recommended \cite{BC3,CB} tool to analyse non-linear operators. To overcome the limitations, we neglect the boundary conditions, extend the operator to an infinite grid and assume that we can linearise the operator locally (we do this by freezing the coefficients). LFA measures the largest amplification factor on high-frequency errors,
for example if there is a smoothing rate of 0.8 this means that the high-frequency
errors are damped by at least 20\%. We initially must derive formulas for the approximation error at each pixel in our 5-point stencil.

{\bf Error forms.} 
Using the definition of $f_{i,j}$, we can rewrite (\ref{eqn:rcabcd})  as
\begin{equation}\label{eqn:abcd}
A_{i,j}\phi_{i+1,j}+B_{i,j}\phi_{i-1,j}+C_{i,j}\phi_{i,j+1}+D_{i,j}\phi_{i,j-1}-S_{i,j}\phi_{i,j}=f_{i,j},
\end{equation}
where we fix $A_{i,j},B_{i,j},C_{i,j}$ and $D_{i,j}$ based on a previous iteration.
The GSLEX I-II and NEWT I-II schemes all work
in a lexicographic manner, and so if we
denote the previous iteration as the $k$-th we can rewrite (\ref{eqn:abcd}) as
\begin{equation}\label{eqn:abcdtimepixel}
A_{i,j}\phi_{i+1,j}^{(k)}+B_{i,j}\phi_{i-1,j}^{(k+1)}+C_{i,j}\phi_{i,j+1}^{(k)}+
D_{i,j}\phi_{i,j-1}^{(k+1)}-S_{i,j}\phi_{i,j}^{(k+1)}=f_{i,j},
\end{equation}
and we obtain the error form by
subtracting (\ref{eqn:abcdtimepixel}) from (\ref{eqn:abcd})
\begin{equation}\label{eqn:errorformpixel}
A_{i,j}e_{i+1,j}^{(k)}+B_{i,j}e_{i-1,j}^{(k+1)}+C_{i,j}e_{i,j+1}^{(k)}+
D_{i,j}e_{i,j-1}^{(k+1)}-S_{i,j}e_{i,j}^{(k+1)}=0,
\end{equation}
Using a similar argument, we obtain the following error form for the line smoothers GSLINE I-II
\begin{equation}\label{eqn:errorform}
A_{i,j}e_{i+1,j}^{(k+1)}+B_{i,j}e_{i-1,j}^{(k+1)}+C_{i,j}e_{i,j+1}^{(k)}+
D_{i,j}e_{i,j-1}^{(k+1)}-S_{i,j}e_{i,j}^{(k+1)}=0,
\end{equation}
where $e_{i,j}^{(k)}=\phi_{i,j}-\phi_{i,j}^{(k)}$ and $e_{i,j}^{(k+1)}=
\phi_{i,j}-\phi_{i,j}^{(k+1)}$.

{\bf Local Fourier Analysis}.\label{sec:LFA}
Define a general Fourier component by
\[
F_{\theta_{1},\theta_{2}}(x_{i},y_{j})=\exp\left( 2\pi {\bf i}\frac{\theta_{1}i}{n} \right)\cdot
\exp\left( 2\pi {\bf i}\frac{\theta_{2}j}{m} \right) = \exp\left( {\bf i}
\frac{\alpha_{1}x_{i}}{h_{x}} \right)\cdot \exp\left( {\bf i}\frac{\alpha_{2}y_{j}}{h_{y}} \right),
\]
where $\alpha_{1}=\frac{2\theta_{1}\pi}{n}$ and $\alpha_{2}=\frac{2\theta_{2}\pi}{m}$
and ${\bf i}$ is the imaginary unit. Note that $\alpha_{1},\alpha_{2}\in[-\pi,\pi]$.
If we assume for simplicity that the image is square and hence $n=m$,
we first expand
\[
e^{(k+1)}_{i,j}=\sum^{n/2}_{\theta_{1},\theta_{2}=-n/2}\psi^{(k+1)}_{\theta_{1},
\theta_{2}}F_{\theta_{1},\theta_{2}}(x_{i},y_{j}),
\hspace{0.5in}
e^{(k)}_{i,j}=\sum^{n/2}_{\theta_{1},\theta_{2}=-n/2}\psi^{(k)}_{\theta_{1},
\theta_{2}}F_{\theta_{1},\theta_{2}}(x_{i},y_{j}),
\]
in Fourier components and define the smoothing rate $\hat{\mu}_{i,j}$ by \cite{BC30,BC13}
\[
\hat{\mu}_{i,j}=\max_{\theta_{1},\theta_{2}}\mu(\theta_{1},\theta_{2})=
\max_{\theta_{1},\theta_{2}}\left|\frac{\psi_{\theta_{1},\theta_{2}}^{(k+1)}}{\psi_{\theta_{1},\theta_{2}}^{(k)}}\right|,
\]
in the high-frequency range where $(\alpha_{1},\alpha_{2})=(\frac{2\theta_{1}\pi}{n},\frac{2\theta_{2}\pi}{n})\in [-\pi,\pi)^{2}\backslash
[-\frac{\pi}{2},\frac{\pi}{2})^{2}$.
Since $\hat{\mu}_{i,j}$ is pixel dependent (non-linear problems), we may also call it the
amplification factor associated with $(i,j)$.

{\bf Smoothing rates.} 
For the GSLEX I-II, NEWT I-II smoothers,
using   (\ref{eqn:errorformpixel}) and (\ref{eqn:errorform}), we obtain error amplification at pixel $(i,j)$
\[
\hat{\mu}_{i,j}=\max_{\theta_{1},\theta_{2}}\mu(\theta_{1},\theta_{2})=\max_{\alpha_{1},\alpha_{2}}\left|\frac{A_{i,j}e^{{\bf i}\alpha_{1}}+C_{i,j}e^{{\bf i}\alpha_{2}}}{B_{i,j}e^{-{\bf i}\alpha_{1}}+D_{i,j}e^{-{\bf i}\alpha_{2}}-S_{i,j}}\right|,
\]
and similarly for the GSLINE I-II smoothers we have
\begin{equation}\label{eqn:gslinelfa}
\hat{\mu}_{i,j}=\max_{\theta_{1},\theta_{2}}\mu(\theta_{1},\theta_{2})=
\max_{\alpha_{1},\alpha_{2}}\left|\frac{C_{i,j}e^{{\bf i}\alpha_{2}}}{A_{i,j}e^{{\bf i}\alpha_{1}}+
B_{i,j}e^{-{\bf i}\alpha_{1}}+D_{i,j}e^{-{\bf i}\alpha_{2}}-S_{i,j}}\right|.
\end{equation}


{\bf Comparison of smoothing rates for all smoothers.}
We consider two different measures of the smoothing rates;
 the maximum and average 
over all pixels $(i,j)$.
We define these in the obvious way as
\vspace{-0.1in}
\[
\tilde{\mu}_{\text{max}}=\max_{i,j}\hat{\mu}_{i,j}=\max_{i,j}\max_{\theta_{1},
\theta_{2}}\mu(\theta_{1},\theta_{2})
\quad \mbox{and}\quad
\tilde{\mu}_{\text{avg}}=\frac{\sum_{i,j}\hat{\mu}_{i,j}}{n^{2}}=
\frac{\sum_{i,j}\max_{\theta_{1},\theta_{2}}\mu(\theta_{1},\theta_{2})}{n^{2}}.
\]
Each of the smoothers was implemented in Algorithm 1 on the image in Figure~\ref{fig:abcdcompare}(a)
with a V-cycle ($\gamma=1$) and using a $1024\times 1024$ resolution image as
the finest grid and a $32\times 32$ image as the coarsest grid and in Table~\ref{tab:compare} we give
$\tilde{\mu}_{\max}$ and $\tilde{\mu}_{\text{avg}}$
for the Rada-Chen and Spencer-Chen models.
\begin{table}[h]
\begin{tabular}{l|c|c|c|r}
\multirow{2}{*}{Smoother } &  \multicolumn{2}{c|}{Rada-Chen} & \multicolumn{2}{c}{Spencer-Chen}\\
& $\tilde{\mu}_{\text{max}}$ & $\tilde{\mu}_{\text{avg}}$
& $\tilde{\mu}_{\text{max}}$ & $\tilde{\mu}_{\text{avg}}$  \\
\hline
\vspace{-0.35cm}&&&\\
GSLINE-I & 0.9997 & 0.4800 & 0.9990 & 0.4586 \\
GSLINE-II & 0.9997 & 0.3782 & 1.0000* & 0.4893 \\
GSLEX-I & 0.9978 & 0.5807 & 0.9927 & 0.5269 \\
GSLEX-II & 1.0000* & 0.5244 & 0.9996 & 0.5512 \\
NEWT-I & 0.9985 & 0.5642 & 0.9595 & 0.4839 \\
NEWT-II & 0.9999 & 0.5749 & 0.9950 & 0.5133 \\
\end{tabular}
\caption{Smoothers and the associated maximum and average smoothing rates for the Rada-Chen and Spencer-Chen models. * due to rounding.}\label{tab:compare}
\end{table}

In the spirit of previous works \cite{BC1}, for any of these smoothers, one would quote
$\tilde{\mu}_{\text{avg}}$, and although this appears to be an excellent rate in all cases,
it is the rate $\tilde{\mu}_{\text{max}}$ that determines the
multigrid convergence \cite{Napov2011}. We therefore choose to focus on $\tilde{\mu}_{\text{max}}$. Table~\ref{tab:compare} shows us that
$\tilde{\mu}_{\max}$ is better for the global smoothers compared to the local smoothers, this is in agreement with the results in \cite{BC1}. However, the maximum smoothing rate of all of the smoothers is bad and so they cannot be implemented
in a successful multigrid scheme. We look to improve the maximum smoothing
rate of one of the better
schemes to obtain a smoother which can be implemented successfully.
In the next section we will see that the problem
is due to discontinuous coefficients in the numerical schemes,
and so we look to \cite{CB1, KC1}
which recommend the use of line smoothers rather than a
pixel-by-pixel update approach. We therefore choose the GSLINE-I smoother and review its performance for the
Rada-Chen model in detail to see if we can improve the maximum smoothing rate of 0.9997.
The same approach will be applied to the Spencer-Chen model and the results
will be quoted at the end of the next section.

{\bf Algorithm 1.} In future discussions, when we compare other algorithms with Algorithm 1, this will be the FAS algorithm with GSLINE-I as smoother.

\section{Non-linear multigrid Algorithm 2}
We now consider how to improve the smoothers above to obtain a smoothing rate which is acceptable.
This leads to our new hybrid smoothers and the resulting multigrid Algorithms 2 and 3.
\subsection{An idea of adaptive iterative schemes}
\label{sec:adaptiveschemes}
\begin{figure}[htb!]
\centering
\includegraphics[scale=0.25]{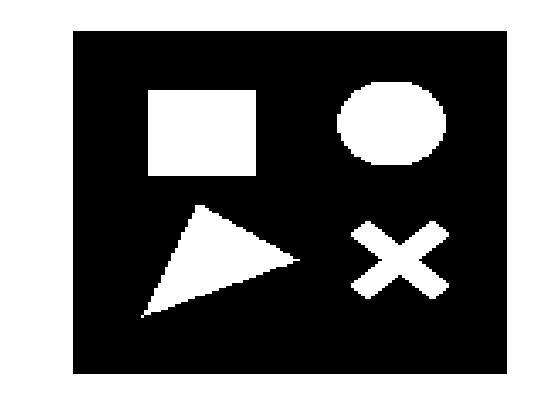}
\includegraphics[scale=0.332]{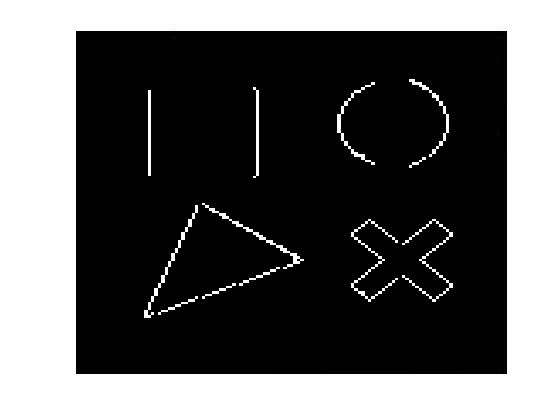}
\includegraphics[width=.325\textwidth]{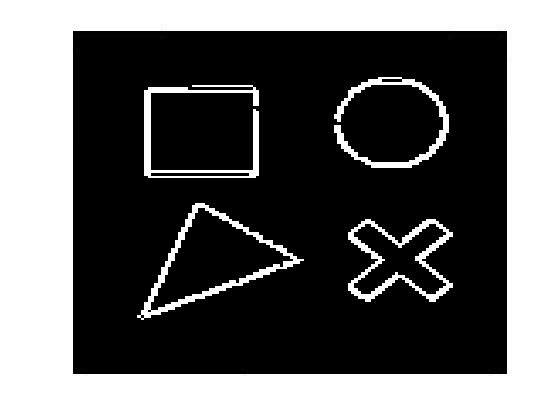}\\
(a)\hspace*{4.5cm}  (b)\hspace*{4.5cm}(c)
\caption{(a) Original image, (b) Pixels with a smoothing rate over 0.6
are indicated in white,
(c) Pixels in white are those where one of  the $A_{i,j},B_{i,j}, C_{i,j}$ or $D_{i,j}$
values differs from the others by a factor of $50\%$ or more.}
\label{fig:abcdcompare}
\end{figure}

To gain more insight into the rates in Table~\ref{tab:compare},
we first look only at those pixels $(i,j)$ which have a large amplification factor.
In Figure~\ref{fig:abcdcompare}(a) we show the original image on
which the rate was measured and in Figure~\ref{fig:abcdcompare}(b) the corresponding binary plot of those pixels
where the amplification factor  is above $0.6$.
We see that the smoother performs poorly at the edges of objects in the image, a phenomenon also observed in \cite{CB} where it was determined that the rate is poor
due to the restriction and interpolation operators performing poorly at these points.

There are two approaches that have been taken to address the
poor smoothing rate at edges; the first is the use of adaptive high
order intergrid transfer operators \cite{CB} and the second is to apply
extra smoothing steps at those edge points \cite{CB5, CB9,CB}.
We prefer the second approach as the intergrid operators perform well
generally and for ease of implementation in the current framework the
second approach is best.
The conventional solution when doing extra smoothing steps would be to simply implement the {\it same} smoother
many more times at those edge pixels to obtain a lower smoothing rate, however we shall develop a {\it different} scheme to be used at these pixels which has an improved smoothing rate.
In any case, we must first   identify those pixels which contribute large amplification factors
without needing to calculate
$\hat{\mu}_{i,j}$ each time, which would be computationally expensive.
In Table~\ref{tab:abcd} we have selected the pixels in the image from
Figure~\ref{fig:abcdcompare}(a) which give 10 of the largest amplification factors
and list the values of $A_{i,j},B_{i,j},C_{i,j}$ and $D_{i,j}$ at these pixels.
\begin{table}[htb!]
\begin{minipage}{0.43\textwidth}
\resizebox{\columnwidth}{!}{%
\begin{tabular}{c|c|c|r|r|r|r}
$i$ & $j$ & $\hat{\mu}_{i,j}$ & \multicolumn{1}{c|}{$A_{i,j}$} & \multicolumn{1}{c|}{$B_{i,j}$} & \multicolumn{1}{c|}{$C_{i,j}$} & \multicolumn{1}{c}{$D_{i,j}$} \\
\hline
\rule{0pt}{2.5ex}
46 & 23 & 0.9997 & 202 & 202 & 137391 & 35 \\
45 & 23 & 0.9995 & 202 & 202 & 77788 & 35 \\
25 & 23 & 0.9931 & 209 & 220 & 5545 & 36 \\
42 & 112 & 0.9889 & 2263 & 1802 & 78959 & 842 \\
44 & 82 & 0.9605 & 20 & 626 & 558 & 22 \\
\end{tabular}%
}
\end{minipage}\qquad 
\begin{minipage}{0.5\textwidth}
\resizebox{\columnwidth}{!}{%
\begin{tabular}{c|c|c|r|r|r|r}
$i$ & $j$ & $\hat{\mu}_{i,j}$ & \multicolumn{1}{c|}{$A_{i,j}$} & \multicolumn{1}{c|}{$B_{i,j}$} & \multicolumn{1}{c|}{$C_{i,j}$} & \multicolumn{1}{c}{$D_{i,j}$} \\
\hline
\rule{0pt}{2.5ex}
44 & 112 & 0.9591 & 79987 & 6659 & 168919 & 6736 \\
97 & 103 & 0.9551 & 3228 & 105968 & 72894 & 3203 \\
80 & 60 & 0.9312 & 7937 & 424357 & 400718 & 27651 \\
73 & 90 & 0.8756 & 29221 & 1426471 & 170469 & 21920 \\
73 & 105 & 0.8750 & 321703 & 24343 & 242663 & 32126 \\
\end{tabular}%
}
\end{minipage}
\caption{The pixels with $10$ of the largest smoothing rates with the corresponding values of $A_{i,j},B_{i,j},C_{i,j}$ and
$D_{i,j}$.}\label{tab:abcd}
\end{table}

\vspace{-0.1in}
A pattern emerges that at these edge pixels (jumps) at least one of the values
of $A_{i,j}, B_{i,j}, C_{i,j}$ and $D_{i,j}$ is significantly different to
the others, Figure~\ref{fig:abcdcompare}(c) shows those pixels
where they differ by 50\% (i.e. $\max (A_{i,j},B_{i,j},C_{i,j},D_{i,j}) / \min (A_{i,j},B_{i,j},C_{i,j},D_{i,j}) > 1.5$).            
%

\begin{defn}
We can identify the edge pixels as those where at least
one of $A_{i,j},B_{i,j}, C_{i,j}$ or $D_{i,j}$ differs significantly from the others,
this is precisely the set of jumps in the coefficients of (\ref{eqn:rcabcd}),
we denote this set by $\mathcal{D}$. For the set of pixels where $A_{i,j},B_{i,j}, C_{i,j}$ or $D_{i,j}$ are relatively
similar we denote it as $\Omega\backslash\mathcal{D}$.
\end{defn}

We compare the maximum and average smoothing rates over $\mathcal{D}$ and
$\Omega\backslash\mathcal{D}$ below:
\begin{equation}
\centering
\fbox{\begin{minipage}{0.56\textwidth}
\begin{tabular}{c|c|c||c|c}
Smoother & $\tilde{\mu}_{\max\mathcal{D}}$ &
$\tilde{\mu}_{\text{avg }\mathcal{D}}$ &
$\tilde{\mu}_{\max\Omega\backslash\mathcal{D}}$ &
$\tilde{\mu}_{\text{avg }\Omega\backslash\mathcal{D}}$ \\
\hline
\vspace{-0.4cm}&&\\
GSLINE-I & 0.9997 & 0.5121 & 0.7705 & 0.4386 \\
\end{tabular}
\end{minipage}}\label{tab:compnonD}
\end{equation}
We see that the maximum amplification factor over $\Omega\backslash\mathcal{D}$ of
0.7705 would mean that the number of
iterations required to reduce the high-frequency errors by 90\% reduces
from 7675 to 9. We now focus on reducing the amplification factor for the pixels of $\mathcal{D}$.

{\bf Classifying the jumps}.
There are $14$ possible cases to consider where one of the
coefficients $A_{i,j},B_{i,j},C_{i,j}$ or $D_{i,j}$ is relatively larger (L) or smaller (S) than
the others, these are all shown below:
\begin{equation}
\resizebox{0.65\columnwidth}{!}{%
\begin{minipage}{0.35\textwidth}
\begin{tabular}{c|c|c|c|c}
Case \# & $A_{i,j}$ & $B_{i,j}$ & $C_{i,j}$ & $D_{i,j}$\\
\hline
\vspace{-0.4cm}&&&\\
1 & S & L & L & S \\
2 & S & L & S & L \\
3 & L & S & L & S \\
4 & L & S & S & L \\
5 & L & L & S & S \\
6 & S & S & L & L \\
7 & L & S & S & S \\
\end{tabular}
\end{minipage}\qquad 
\begin{minipage}{0.35\textwidth}
\begin{tabular}{c|c|c|c|c}
Case \# & $A_{i,j}$ & $B_{i,j}$ & $C_{i,j}$ & $D_{i,j}$\\
\hline
\vspace{-0.4cm}&&&\\
8 & S & S & L & S \\
9 & S & L & S & S \\
10 & S & S & S & L \\
11 & L & L & S & L \\
12 & L & S & L & L \\
13 & L & L & L & S \\
14 & S & L & L & L \\
\end{tabular}
\end{minipage}
}
\label{tab:compareabcd}
\end{equation}
We can now label each pixel in $\mathcal{D}$ as one of the cases from $1$ to $14$. The choice of label $L$ or $S$ for a coefficient will be dependent on the coefficients at each pixel. Typically, if the largest coefficient is 50\% larger than the smallest we group the coefficients as large or small by K-means or some other classification method. For a pixel in $\mathcal{D}$, we now look to adapt the iterative
scheme (\ref{eqn:abcd}) for each of these cases to give a scheme which has a
better smoothing rate than implementing GSLINE-I directly.
In the interests of brevity, we consider Case $1$ in detail and will generalise the results
to  other cases next.

\subsubsection{An adapted iterative scheme and its LFA form}
Our aim is to propose a new iteration scheme
which leads to a smaller smoothing rate by  the LFA.
For Case 1 pixels, $A_{i,j}$ and $D_{i,j}$ are relatively small and $B_{i,j}$ and
$C_{i,j}$ are relatively large.  We can rewrite (\ref{eqn:abcd}) as
\[
B_{i,j}\phi_{i-1,j}+C_{i,j}\phi_{i,j+1}-S_{i,j}\phi_{i,j}=f_{i,j}-A_{i,j}\phi_{i+1,j}-
D_{i,j}\phi_{i,j-1},
\]
by moving the small terms to the right hand side. We now look to solve
$\phi_{i-1,j},\phi_{i,j+1}$ and $\phi_{i,j}$ as a coupled system.
We can rewrite this scheme, with the iteration number indicated, as
\begin{equation}\label{eqn:genlfacase12}
B_{i,j}\phi_{i-1,j}^{(k+1)}+C_{i,j}\phi_{i,j+1}^{(k+1)}-S_{i,j}\phi_{i,j}^{(k+1)}=
f_{i,j}- A_{i,j}\phi_{i+1,j}^{(k)}-D_{i,j}\phi_{i,j-1}^{(k)}.
\end{equation}

The amplification factor for such a scheme is
\begin{equation}\label{eqn:smratecase1}
\hat{\mu}_{i,j}=\max_{\theta_{1},\theta_{2}}\mu(\theta_{1},\theta_{2})=
\max_{\alpha_{1},\alpha_{2}}\frac{|A_{i,j}e^{{\bf i}\alpha_{1}}+D_{i,j}e^{-{\bf i}
\alpha_{2}}|}{|S_{i,j}-B_{i,j}e^{-{\bf i}\alpha_{1}}-C_{i,j}e^{{\bf i}\alpha_{2}}|},
\end{equation}
derived as in \S\ref{sec:LFA}. In fact,
we see the following improvements to the maximum and average smoothing rates for all of the Case 1 pixels by using the adapted iterative
scheme (\ref{eqn:genlfacase12}) rather than the GSLINE-I smoother in (\ref{eqn:system})
\[
 \tilde{\mu}_{\max} =0.9863, \
 \tilde{\mu}_{\text{avg}}=0.7174\quad
 \Longrightarrow\quad
 \tilde{\mu}_{\max} =0.7324, \
 \tilde{\mu}_{\text{avg}}= 	0.3013
\]
Reducing the smoothing rate from 0.9863 to 0.7324 is dramatic; exemplified by the fact that to reduce high-frequency errors
by 90\% for Case 1 pixels with GSLINE-I we
would have required 167 iterations but now we need just 8.
Hence, now we know that the scheme (\ref{eqn:genlfacase12}) gives us a
better smoothing rate than GSLINE-I at these pixels.

\vspace{-0.1in}
\subsubsection{Adapted schemes for all cases
of (\ref{tab:compareabcd}) and their rates by LFA}
Using the central idea of lagging the small terms
in (\ref{tab:compareabcd}) (between $1$ and $3$ terms),
we can derive adapted schemes for all cases 
in the same manner as for Case 1
previously. In Table~\ref{tab:comparerates2} we give the comparison of the
maximum smoothing rate of GSLINE-I, $\mu_{GSLINE}$, with the maximum
smoothing rate of the adapted schemes $\mu_{adapted_{1}}$.
\begin{table}[htb!]
\resizebox{0.7\columnwidth}{!}{
\begin{minipage}{0.35\textwidth}
\begin{tabular}{c|c|c|c}
Case \#& $\mu_{GSLINE}$ & $\mu_{adapted_1}$ & \\
\hline
\vspace{-0.4cm}&&\\
1 & 0.9863 & 0.7324 & $\diamondsuit$ \\
2 & 0.6259 & 0.8515 & $\diamondsuit$ \\
3 & 0.9900 & 0.7418 & $\diamondsuit$ \\
4 & 0.6408 & 0.7415 & $\diamondsuit$ \\
5 & 0.7105 & 1.0000 & $\Box$ \\
6 & 0.9524 & 1.0000 & $\Box$ \\
7 & 0.9592 & 0.9536 & $\diamondsuit$ \\
\end{tabular}
\end{minipage}\qquad 
\begin{minipage}{0.35\textwidth}
\begin{tabular}{c|c|c|c}
Case \#& $\mu_{GSLINE}$ & $\mu_{adapted_1}$  \\
\hline
\vspace{-0.4cm}&&\\
8 & 0.9997 & 0.9569 & $\diamondsuit$ \\
9 & 0.9481 & 0.9426 & $\diamondsuit$ \\
10 & 0.8935 & 0.9640 & $\diamondsuit$ \\
11 & 0.2693 & 0.2693 & $\spadesuit$ \\
12 & 0.7729 & 0.2663 & $\spadesuit$ \\
13 & 0.9865 & 0.2704 & $\spadesuit$ \\
14 & 0.5993 & 0.2706 & $\spadesuit$
\end{tabular}
\end{minipage}
}
\caption{Comparison of the maximum amplification factors using GSLINE-I and the
adapted iterative schemes for each case. The $\Box$-cases are the decoupled cases which give a rate of precisely 1, as remarked, the $\diamondsuit$-cases have minor or no improvement in the smoothing rate and the $\spadesuit$-cases have a good final rate.}\label{tab:comparerates2}
\end{table}

The results from Table \ref{tab:comparerates2}
fall into $3$ categories:
\begin{itemize}
\item[] $\spadesuit$-cases, where only one term is lagged and the improvements are remarkable. 
    This gives a promising indication that the lagging of particular terms in certain cases can improve the smoothing rate.
     This motivates our next step.
  
  \item[] 
  $\diamondsuit$-cases, where either $2$ or $3$ terms are lagged.
  We see either only a minor improvement to an already
 high rate or the rate has actually worsened. 

 \item[] $\Box$-cases, where $2$ terms are lagged and we see the worst results: a smoothing rate of 1.0000 is attained 
 for cases $5,6$ in Table~\ref{tab:comparerates2}.
Below we prove analytically that for Case $6$ pixels the smoothing rate when using the adapted scheme will always be precisely 1.

Case $6$ pixels have the LFA form $\hat{\mu}_{i,j}=
\max_{\alpha_{1},\alpha_{2}}\frac{|A_{i,j}e^{{\bf i}\alpha_{1}}+B_{i,j}e^{-{\bf i}
\alpha_{1}}|}{|S_{i,j}-C_{i,j}e^{{\bf i}\alpha_{2}}-D_{i,j}e^{-{\bf i}\alpha_{2}}|},$
and we see a decoupling in the maximisation with respect to $\alpha_{1}$ and $\alpha_{2}$ which allows us to rewrite this as

\resizebox{\columnwidth}{!}{$
\begin{gathered}
\begin{aligned}
\hat{\mu}_{i,j}&
=\frac{\max\limits_{\alpha_{1}}\Big|A_{i,j}e^{{\bf i}\alpha_{1}}+B_{i,j}e^{-{\bf i}
\alpha_{1}}\Big|}{\min\limits_{\alpha_{2}}\Big|S_{i,j}-C_{i,j}e^{{\bf i}\alpha_{2}}-D_{i,j}e^{-{\bf i}\alpha_{2}}\Big|}= \frac{\max\limits_{\alpha_{1}}\Big|(A_{i,j}+B_{i,j})\cos(\alpha_{1})+{\bf i}(A_{i,j}-B_{i,j})\sin(\alpha_{1})\Big|}{\min\limits_{\alpha_{2}}\Big|\left[A_{i,j}+B_{i,j}+C_{i,j}(1-\cos(\alpha_{2}))+D_{i,j}(1-\cos(\alpha_{2}))\right]+{\bf i}(C_{i,j}-D_{i,j})\sin(\alpha_{2})\Big|} \\
\rule{0pt}{7ex}
&= \frac{\sqrt{\max\limits_{\alpha_{1}}\Big[A_{i,j}^{2}+B_{i,j}^{2}+
2A_{i,j}B_{i,j}cos(2\alpha_{1})\Big]}}{\sqrt{\min\limits_{\alpha_{2}}\Big[\left[A_{i,j}+B_{i,j}+C_{i,j}(1-\cos(\alpha_{2}))+D_{i,j}(1-\cos(\alpha_{2}))\right]^{2}+(C_{i,j}-D_{i,j})^{2}\sin(\alpha_{2})^{2}\Big]}}
\quad = \quad\frac{(A_{i,j}+B_{i,j})^{2}}{(A_{i,j}+B_{i,j})^{2}} = 1,\\
\end{aligned}
\end{gathered}$}
attained at $(\alpha_{1},\alpha_{2}) = (-\pi,0)\in[-\pi,\pi)^{2}\backslash
[-\frac{\pi}{2},\frac{\pi}{2})^{2}$. Similarly we have $\hat{\mu}_{i,j}=1$ for Case $5$ too.

We claim that
it is necessary to have both of $\alpha_{1}$ and $\alpha_{2}$ in the numerator or
denominator of the LFA formulation to ensure a low smoothing rate.
We note that for
 Cases $5$ and $6$ this is not the case.
 \end{itemize}

We now focus on improving the $\diamondsuit$-cases and  the Case 8 in particular and
its LFA to motivate us on how to proceed i.e.
to see whether an alternative adaptation to the iterative scheme
gives a better smoothing rate. The results apply to $\Box$-cases also.

\begin{figure}[htb!]
\mbox{}\hspace*{-0.03cm}
\includegraphics[height=5.5cm,width=\textwidth]{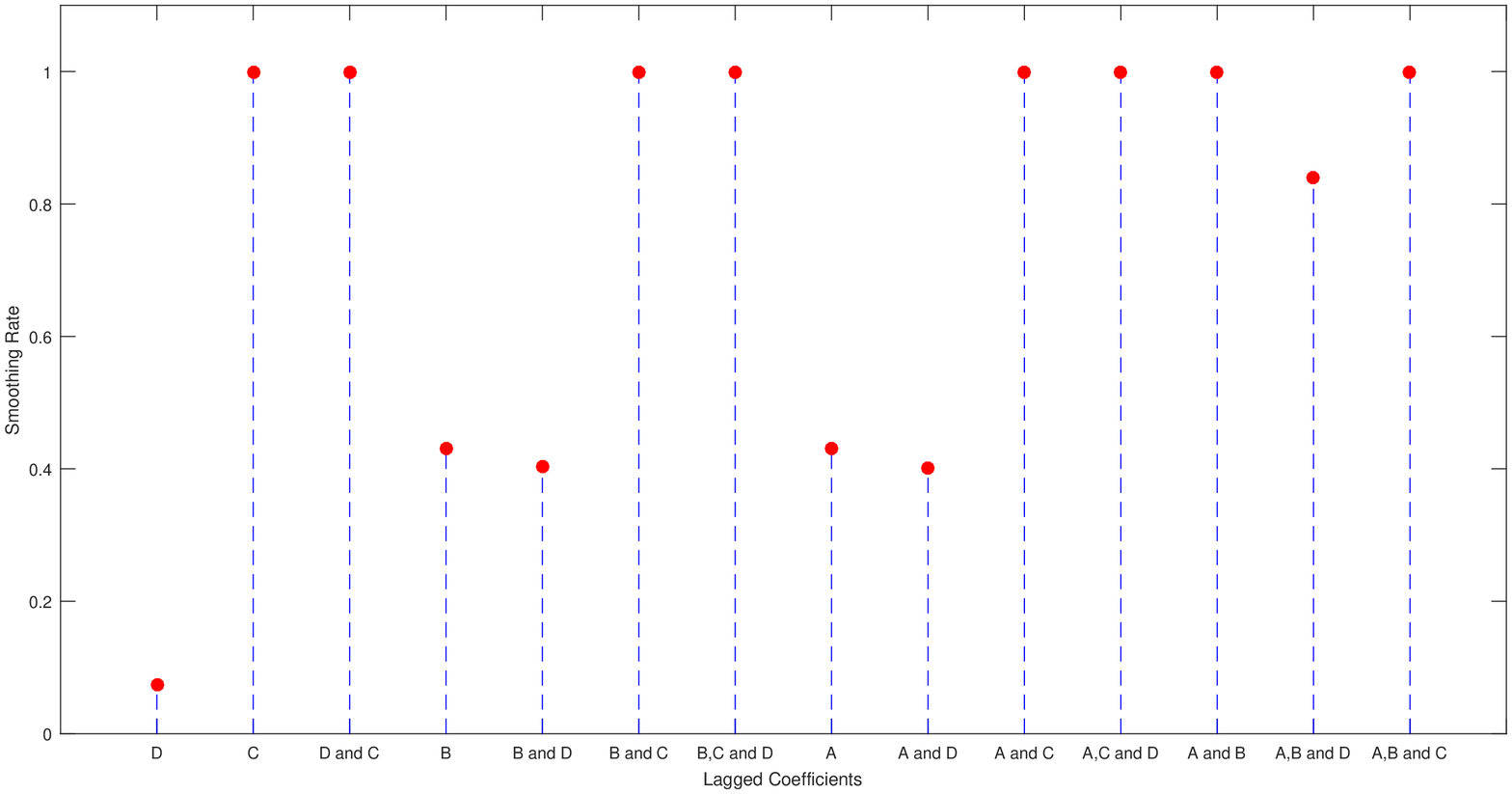}
\vspace{-0.3in}
\caption{Comparison of the smoothing rate for the Case 8 pixel with the worst smoothing rate when different coefficient terms are lagged. In this case, $A_{i,j} = 202, B_{i,j} = 202, C_{i,j} = 137391$ and $D_{i,j} = 35$ (Table~\ref{tab:abcd}).}\label{fig:lfacompare}
\end{figure}

{\bf Improving the adapted scheme for Case 8}. A pixel which is labelled as Case $8$ is one where $A_{i,j}, B_{i,j}, D_{i,j}$ are relatively small and $C_{i,j}$ is relatively large. Using the previous method we would devise a scheme where the terms with coefficients $A_{i,j}, B_{i,j}, D_{i,j}$ would be lagged at time step $k$ and the term with coefficient $C_{i,j}$ would be updated to time step $k+1$. We pick the particular Case 8 pixel from Table~\ref{tab:abcd} which has the worst smoothing rate and in Figure~\ref{fig:lfacompare} we look at the smoothing rate for the scheme (\ref{eqn:abcd}) with
different coefficients lagged.

This shows that the best rate is achieved when just the smallest of the coefficients ($D_{i,j}$) is lagged. Even the lagging of two of the smallest coefficients gives an improvement on lagging all three. This gives some indication that the smoothing rate is best when the smallest coefficient is lagged and this has proven to be the case in every one of the many examples which the authors have tried. It would be an interesting piece of future work to prove that this must be true analytically.

Hence we propose to lag just the smallest of the coefficients in a modified scheme for all cases.

\subsubsection{Improved adapted schemes for all cases}
We re-consider the $\Diamond$ and $\Box$-cases which have more than one relatively small coefficient.
Lagging only the smallest coefficient, the LFA forms simplify to those of
Cases 11--14 and we expect major improvements. In Table~\ref{tab:comparerates3} we compare the maximum smoothing
rate of GSLINE-I, $\mu_{GSLINE}$, for these cases with the maximum smoothing rate
of an improved, adapted iterative scheme which lags only the smallest coefficient
$\mu_{adapted_{2}}$.
\begin{table}[htb!]
\resizebox{0.7\columnwidth}{!}{
\begin{minipage}{0.35\textwidth}
\begin{tabular}{c|c|c}
Case \#& $\mu_{GSLINE}$ & $\mu_{adapted_2}$  \\
\hline
\vspace{-0.4cm}&&\\
1 & 0.9863 & 0.4467 \\
2 & 0.6259 & 0.4398 \\
3 & 0.9900 & 0.4280 \\
4 & 0.6408 & 0.4468 \\
5 & 0.7105 & 0.4659 \\
6 & 0.9524 & 0.4547 \\
7 & 0.9592 & 0.4789 \\
\end{tabular}
\end{minipage}\qquad 
\begin{minipage}{0.35\textwidth}
\begin{tabular}{c|c|c}
Case \#& $\mu_{GSLINE}$ & $\mu_{adapted_2}$  \\
\hline
\vspace{-0.4cm}&&\\
8 & 0.9997 & 0.4779 \\
9 & 0.9481 & 0.4716 \\
10 & 0.8935 & 0.4749 \\
11 & 0.2693 & 0.2693 \\
12 & 0.7729 & 0.2663 \\
13 & 0.9865 & 0.2704 \\
14 & 0.5993 & 0.2706 \\
\end{tabular}
\end{minipage}
}
\caption{Comparison of the maximum amplification factors using GSLINE-I and the adapted iterative schemes for each case with just the smallest coefficient term lagged.}\label{tab:comparerates3}
\end{table}

As expected, 
there is a significant improvement in the smoothing rate in all cases when we lag
 just the smallest coefficient, it also makes implementation faster as we now
 consider just 4 cases of possible lagged coefficients rather than 14 and therefore have
 only 4 iterative schemes to consider. Taking our guidance from these results, we propose two hybrid smoothers which both perform standard smoothing iterations on pixels of $\Omega\backslash\mathcal{D}$ and perform non-standard adapted iterative schemes on the pixels in $\mathcal{D}$. 
 
Based on the above pixel-wise motivating tests, we now present two 
iterative strategies for our new smoothers.
The first smoother is natural: for each pixel $(i,j)$, in $\mathcal{D}$, all of the directly connected neighbouring pixels are collectively updated except the term with the smallest coefficient.
That is, Hybrid Smoother 1 uses block structure Vanka-type smoothing schemes \cite{Rodrigo:16,Vanka:86} to update the pixels in $\mathcal{D}$.
The potential drawback is that previously updated pixels may enter to the next group of (potentially multiple) updates, making subsequent analysis intractable.
Hence our second smoother, 
 denoted by `Hybrid Smoother 2',  incorporates partial line smoothing operations at pixels in $\mathcal{D}$ and only pixels that are the same line as $(i,j)$ are updated. This line by line approach facilitates subsequent analysis.

\subsection{Hybrid Smoother 1} 

Our first hybrid  smoother
updates blocks of pixels at each update, these blocks may overlap. 
This is an overlapping block smoother of Vanka-type \cite{Rodrigo:16,Vanka:86}.
Once again we start with the set $\mathcal{D}$ of pixels with jumping coefficients. 
 For brevity, we will detail the derivation of the iterative scheme for pixels in $\mathcal{D}$ for which $A_{i,j}$ is smallest. We will then state the schemes for the other laggings (derived in the same manner).
 
{\bf $A_{i,j}$ lagged.} 
The lagging of coefficient $A_{i,j}$ in equation (\ref{eqn:abcd}) gives rise to the iterative scheme
\vspace{-0.1in}
\begin{equation}\label{eqn:abcd2}
A_{i,j}\phi_{i+1,j}^{(k)}+B_{i,j}\phi_{i-1,j}^{(k+1)}+C_{i,j}\phi_{i,j+1}^{(k+1)}+D_{i,j}\phi_{i,j-1}^{(k+1)}-S_{i,j}\phi_{i,j}^{(k+1)}=f_{i,j},
\end{equation}
We are solving for $\phi_{i-1,j},\phi_{i,j+1},\phi_{i,j-1}$ and
$\phi_{i,j}$ simultaneously and as we have only one equation, we need three more.
We get these by considering (\ref{eqn:abcd}) at the pixels $(i-1,j)$ and
$(i,j+1)$ and $(i,j-1)$, which gives us the three equations
\begin{equation*}
\begin{gathered}
B_{i,j}\phi_{i,j}-S_{i-1,j}\phi_{i-1,j}=f_{i-1,j}-B_{i-1,j}\phi_{i-2,j}-
C_{i-1,j}\phi_{i-1,j+1}-D_{i-1,j}\phi_{i-1,j-1},\\
C_{i,j}\phi_{i,j}-S_{i,j+1}\phi_{i,j+1}=f_{i,j+1}-A_{i,j+1}\phi_{i+1,j+1}-
B_{i,j+1}\phi_{i-1,j+1}-C_{i,j+1}\phi_{i,j+2},\\
D_{i,j}\phi_{i,j}-S_{i,j-1}\phi_{i,j-1}=f_{i,j-1}-A_{i,j-1}\phi_{i+1,j-1}-
B_{i,j-1}\phi_{i-1,j-1}-D_{i,j-1}\phi_{i,j-2},\\
\end{gathered}
\end{equation*}
which have been rearranged to have the $\phi_{i-1,j},\phi_{i,j+1},\phi_{i,j-1}$ and $\phi_{i,j}$
terms on the left hand side. So, using these along with (\ref{eqn:abcd2})
we obtain the system (\ref{eqn:caseAlagged}).

{\bf Scheme with $A_{i,j}$ lagged}:
\begin{equation}\resizebox{.9\hsize}{!}{$\left(\begin{array}{cccc}
-S_{i,j} & B_{i,j} &  C_{i,j} & D_{i,j} \\
B_{i,j} & -S_{i-1,j} & 0& 0 \\
C_{i,j} & 0 & -S_{i,j+1}& 0 \\
D_{i,j} & 0 & 0 & -S_{i,j-1}\\
\end{array}\right)\cdot\left(
\begin{array}{c}
\phi_{i,j}\\
\phi_{i-1,j}\\
\phi_{i,j+1}\\
\phi_{i,j-1}\\
\end{array}\right)
=\left(
\begin{array}{c}
f_{i,j}-A_{i,j}\phi_{i+1,j}\\
f_{i-1,j}-C_{i-1,j}\phi_{i-1,j+1}-D_{i-1,j}\phi_{i-1,j-1}-B_{i-1,j}\phi_{i-2,j}\\
f_{i,j+1}-A_{i,j+1}\phi_{i+1,j+1}-B_{i,j+1}\phi_{i-1,j+1}-C_{i,j+1}\phi_{i,j+2}\\
f_{i,j-1}-A_{i,j-1}\phi_{i+1,j-1}-B_{i,j-1}\phi_{i-1,j-1}-D_{i,j-1}\phi_{i,j-2}\\
\end{array}
\right)$}.\label{eqn:caseAlagged}\end{equation}
%
This system is strictly diagonally dominant and follows the guidance in \cite{BC30} that collective update schemes are better for jumping coefficients. This system also has an arrow structure in the matrix and can be solved very quickly (in 24 operations).

\subsubsection{The adapted iterative schemes for other cases}

Below are the adapted iterative schemes for the cases when $B_{i,j}, C_{i,j}$ or $D_{i,j}$ are lagged, derived in the same manner as previously when $A_{i,j}$ was lagged.

{\bf Scheme with $B_{i,j}$ lagged}:
 \begin{equation}\resizebox{.9\hsize}{!}{$\left(\begin{array}{cccc}
-S_{i,j} & A_{i,j} &  C_{i,j} & D_{i,j} \\
A_{i,j} & -S_{i+1,j} & 0& 0 \\
C_{i,j} & 0 & -S_{i,j+1}& 0 \\
D_{i,j} & 0 & 0 & -S_{i,j-1}\\
\end{array}\right)\cdot\left(
\begin{array}{c}
\phi_{i,j}\\
\phi_{i+1,j}\\
\phi_{i,j+1}\\
\phi_{i,j-1}\\
\end{array}\right)
=\left(
\begin{array}{c}
f_{i,j}-B_{i,j}\phi_{i-1,j}\\
f_{i+1,j}-C_{i+1,j}\phi_{i+1,j+1}-D_{i+1,j}\phi_{i+1,j-1}-A_{i+1,j}\phi_{i+2,j}\\
f_{i,j+1}-A_{i,j+1}\phi_{i+1,j+1}-B_{i,j+1}\phi_{i-1,j+1}-C_{i,j+1}\phi_{i,j+2}\\
f_{i,j-1}-A_{i,j-1}\phi_{i+1,j-1}-B_{i,j-1}\phi_{i-1,j-1}-D_{i,j-1}\phi_{i,j-2}\\
\end{array}
\right)$}.\label{eqn:case12}\end{equation}

{\bf Scheme with $C_{i,j}$ lagged}:
\begin{equation}\resizebox{.9\hsize}{!}{$\left(\begin{array}{cccc}
-S_{i,j} & A_{i,j} &  B_{i,j} & D_{i,j} \\
A_{i,j} & -S_{i+1,j} & 0& 0 \\
B_{i,j} & 0 & -S_{i-1,j}& 0 \\
D_{i,j} & 0 & 0 & -S_{i,j-1}\\
\end{array}\right)\cdot\left(
\begin{array}{c}
\phi_{i,j}\\
\phi_{i+1,j}\\
\phi_{i-1,j}\\
\phi_{i,j-1}\\
\end{array}\right)
=\left(
\begin{array}{c}
f_{i,j}-C_{i,j}\phi_{i,j+1}\\
f_{i+1,j}-C_{i+1,j}\phi_{i+1,j+1}-D_{i+1,j}\phi_{i+1,j-1}-A_{i+1,j}\phi_{i+2,j}\\
f_{i-1,j}-C_{i-1,j}\phi_{i-1,j+1}-D_{i-1,j}\phi_{i-1,j-1}-B_{i-1,j}\phi_{i-2,j}\\
f_{i,j-1}-A_{i,j-1}\phi_{i+1,j-1}-B_{i,j-1}\phi_{i-1,j-1}-D_{i,j-1}\phi_{i,j-2}\\
\end{array}
\right)$}.\label{eqn:case11}\end{equation}

{\bf Scheme with $D_{i,j}$ lagged}:
\begin{equation}\resizebox{.9\hsize}{!}{$\left(\begin{array}{cccc}
-S_{i,j} & A_{i,j} &  B_{i,j} & C_{i,j} \\
A_{i,j} & -S_{i+1,j} & 0& 0 \\
B_{i,j} & 0 & -S_{i-1,j}& 0 \\
C_{i,j} & 0 & 0 & -S_{i,j+1}\\
\end{array}\right)\cdot\left(
\begin{array}{c}
\phi_{i,j}\\
\phi_{i+1,j}\\
\phi_{i-1,j}\\
\phi_{i,j+1}\\
\end{array}\right)
=\left(
\begin{array}{c}
f_{i,j}-D_{i,j}\phi_{i,j-1}\\
f_{i+1,j}-C_{i+1,j}\phi_{i+1,j+1}-D_{i+1,j}\phi_{i+1,j-1}-A_{i+1,j}\phi_{i+2,j}\\
f_{i-1,j}-C_{i-1,j}\phi_{i-1,j+1}-D_{i-1,j}\phi_{i-1,j-1}-B_{i-1,j}\phi_{i-2,j}\\
f_{i,j+1}-A_{i,j+1}\phi_{i+1,j+1}-B_{i,j+1}\phi_{i-1,j+1}-C_{i,j+1}\phi_{i,j+2}\\
\end{array}
\right)$}.\label{eqn:caseDlagged}\end{equation}

\subsubsection{Implementing Hybrid Smoother 1}

To minimise grid sweeps and ensure that all pixels are covered, we use the following pseudo-algorithm for Hybrid Smoother 2:
\begin{enumerate}
\item[{\bf I\ \ }]{Perform GSLINE-I on all lines in the image.}
\item[{\bf II\ }]{For each pixel in $\mathcal{D}$, perform the appropriate scheme of (\ref{eqn:caseAlagged})--(\ref{eqn:caseDlagged}).}
\end{enumerate}
We justify the choice of GSLINE-I in step {\bf I} as it is the recommended smoothing scheme for a problem with jump coefficients \cite{BC30}. Note that the schemes in {\bf II} can overlap the same pixels several times due to the collective updates.

{\bf Algorithm 2.} In future discussion, when we use the Hybrid Smoother 1 in the Full Approximation Scheme, we will call this Algorithm 2.


\subsection{Hybrid Smoother 2}
\label{sec:hyb1}

Our second hybrid smoother first groups pixels in $\mathcal{D}$ 
by whether $A_{i,j}$, $B_{i,j}$, $C_{i,j}$ or $D_{i,j}$ are the smallest and then by the line they are on.
We then perform partial line updates on these groups for $A_{i,j}$, $B_{i,j}$, $C_{i,j}$ or $D_{i,j}$ in sequence along with individual pixel updates on the other pixels, this avoids the overlap encountered in Hybrid Smoother 1.
We note that for pixels in $\Omega\backslash\mathcal{D}$ the LFA tells us that the smoothing rate is acceptable (maximum 0.7705) and therefore we design a smoother which performs cheap GSLEX-I iterations at the pixels of $\Omega\backslash\mathcal{D}$ and performs the lagged scheme on the other pixels. We focus initially on how we propose implementing this for the pixels in $\mathcal{D}$ with $A_{i,j}$ lagged and then we generalise the idea to the laggings of $B_{i,j},C_{i,j}$ and $D_{i,j}$.

{\bf Scheme with $A_{i,j}$ lagged}.
Suppose we focus on a pixel $(i,j) \in\mathcal{D}$ which has coefficient $A_{i,j}$ the smallest. If we lag the $A_{i,j}$ the smoothing rate at this pixel is
\[
\hat{\mu}_{i,j}=
\max_{(\alpha_{1},\alpha_{2})\in[-\pi,\pi)^{2}\backslash
[-\frac{\pi}{2},\frac{\pi}{2})^{2}}\left|\frac{A_{i,j}e^{{\bf i}\alpha_{1}}}{B_{i,j}e^{-{\bf i}\alpha_{1}}+C_{i,j}e^{{\bf i}\alpha_{2}}+
D_{i,j}e^{-{\bf i}\alpha_{2}}-S_{i,j}}\right|
\]
which is precisely the smoothing rate for a line smoother updating from the top row to the bottom row.
In the majority of cases, if pixel $A_{i,j}$ is the smallest, we find that many adjacent pixels on that line also have $A_{i,\cdot}$ the smallest. So we can perform a partial line smoothing on these pixels.

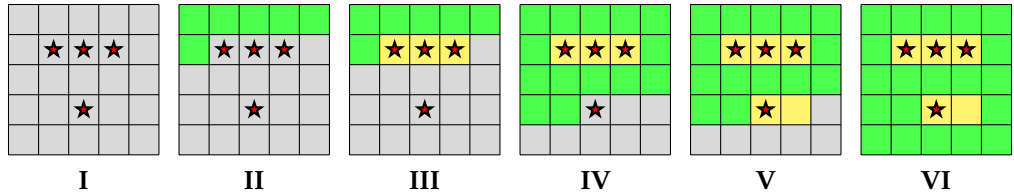
\begin{figure}[htb!]
  \begin{minipage}{0.15\linewidth}
    \centering
      \begin{tikzpicture}[darkstyle/.style={circle,draw,fill=black,minimum size=1,scale = 0.5},scale = 0.4]
      \foreach \x in {0,...,2}
   \foreach \y in {0,...,2}
{
\fill [gray!30] (0, 0 ) rectangle (5,5);
\ngram{0.25}{5}{54}{thick,fill=red}{1}{3};
\ngram{0.25}{5}{54}{thick,fill=red}{2}{3};
\ngram{0.25}{5}{54}{thick,fill=red}{3}{3};
\ngram{0.25}{5}{54}{thick,fill=red}{2}{1};
}
  \foreach \x in {0,...,5}
      \draw (\x,0)--(\x,5) (0,\x)--(5,\x) ;
\end{tikzpicture}

      {\bf {I}}
  \end{minipage}
      \begin{minipage}{0.15\linewidth}
    \centering
  \begin{tikzpicture}[darkstyle/.style={circle,draw,fill=black,minimum size=1,scale = 0.5},scale = 0.4]
      \foreach \x in {0,...,2}
   \foreach \y in {0,...,2}
{
\fill [gray!30] (0, 0 ) rectangle (5,5);
\ngram{0.25}{5}{54}{thick,fill=red}{1}{3};
\ngram{0.25}{5}{54}{thick,fill=red}{2}{3};
\ngram{0.25}{5}{54}{thick,fill=red}{3}{3};
\ngram{0.25}{5}{54}{thick,fill=red}{2}{1};
\fill [green!70] (5 , 5 ) rectangle (0,4);
\fill [green!70] (1 , 4 ) rectangle (0,3);
}
  \foreach \x in {0,...,5}
      \draw (\x,0)--(\x,5) (0,\x)--(5,\x) ;
\end{tikzpicture}
  {\bf {II}}
     \end{minipage}
         \begin{minipage}{0.15\linewidth}
    \centering
    \begin{tikzpicture}[darkstyle/.style={circle,draw,fill=black,minimum size=1,scale = 0.5},scale = 0.4]
      \foreach \x in {0,...,2}
   \foreach \y in {0,...,2}
{
\fill [gray!30] (0, 0 ) rectangle (5,5);
\fill [yellow!70] (1 , 4 ) rectangle (4,3);
\ngram{0.25}{5}{54}{thick,fill=red}{1}{3};
\ngram{0.25}{5}{54}{thick,fill=red}{2}{3};
\ngram{0.25}{5}{54}{thick,fill=red}{3}{3};
\ngram{0.25}{5}{54}{thick,fill=red}{2}{1};
\fill [green!70] (5 , 5 ) rectangle (0,4);
\fill [green!70] (1 , 4 ) rectangle (0,3);
}
  \foreach \x in {0,...,5}
      \draw (\x,0)--(\x,5) (0,\x)--(5,\x) ;
\end{tikzpicture}
      {\bf {III}}
     \end{minipage}
         \begin{minipage}{0.15\linewidth}
    \centering
    \begin{tikzpicture}[darkstyle/.style={circle,draw,fill=black,minimum size=1,scale = 0.5},scale = 0.4]
      \foreach \x in {0,...,2}
   \foreach \y in {0,...,2}
{
\fill [gray!30] (0, 0 ) rectangle (5,5);
\fill [yellow!70] (1 , 4 ) rectangle (4,3);
\ngram{0.25}{5}{54}{thick,fill=red}{1}{3};
\ngram{0.25}{5}{54}{thick,fill=red}{2}{3};
\ngram{0.25}{5}{54}{thick,fill=red}{3}{3};
\ngram{0.25}{5}{54}{thick,fill=red}{2}{1};
\fill [green!70] (5 , 5 ) rectangle (0,4);
\fill [green!70] (1 , 4 ) rectangle (0,3);
\fill [green!70] (4 , 3 ) rectangle (5,4);
\fill [green!70] (5 , 3 ) rectangle (0,2);
\fill [green!70] (2 , 2) rectangle (0,1);
}
  \foreach \x in {0,...,5}
      \draw (\x,0)--(\x,5) (0,\x)--(5,\x) ;
\end{tikzpicture}
      {\bf {IV}}
     \end{minipage}
         \begin{minipage}{0.15\linewidth}
    \centering
    \begin{tikzpicture}[darkstyle/.style={circle,draw,fill=black,minimum size=1,scale = 0.5},scale = 0.4]
      \foreach \x in {0,...,2}
   \foreach \y in {0,...,2}
{
\fill [gray!30] (0, 0 ) rectangle (5,5);
\fill [yellow!70] (1 , 4 ) rectangle (4,3);
\fill [yellow!70]  (2 , 2) rectangle (4,1);
\ngram{0.25}{5}{54}{thick,fill=red}{1}{3};
\ngram{0.25}{5}{54}{thick,fill=red}{2}{3};
\ngram{0.25}{5}{54}{thick,fill=red}{3}{3};
\ngram{0.25}{5}{54}{thick,fill=red}{2}{1};
\fill [green!70] (5 , 5 ) rectangle (0,4);
\fill [green!70] (1 , 4 ) rectangle (0,3);
\fill [green!70] (4 , 3 ) rectangle (5,4);
\fill [green!70] (5 , 3 ) rectangle (0,2);
\fill [green!70] (2 , 2) rectangle (0,1);
}
  \foreach \x in {0,...,5}
      \draw (\x,0)--(\x,5) (0,\x)--(5,\x) ;
\end{tikzpicture}
      {\bf {V}}
     \end{minipage}
         \begin{minipage}{0.15\linewidth}
    \centering
    \begin{tikzpicture}[darkstyle/.style={circle,draw,fill=black,minimum size=1,scale = 0.5},scale = 0.4]
      \foreach \x in {0,...,2}
   \foreach \y in {0,...,2}
{
\fill [gray!30] (0, 0 ) rectangle (5,5);
\fill [yellow!70] (1 , 4 ) rectangle (4,3);
\fill [yellow!70]  (2 , 2) rectangle (4,1);
\ngram{0.25}{5}{54}{thick,fill=red}{1}{3};
\ngram{0.25}{5}{54}{thick,fill=red}{2}{3};
\ngram{0.25}{5}{54}{thick,fill=red}{3}{3};
\ngram{0.25}{5}{54}{thick,fill=red}{2}{1};
\fill [green!70] (5 , 5 ) rectangle (0,4);
\fill [green!70] (1 , 4 ) rectangle (0,3);
\fill [green!70] (4 , 3 ) rectangle (5,4);
\fill [green!70] (5 , 3 ) rectangle (0,2);
\fill [green!70] (2 , 2) rectangle (0,1);
\fill [green!70] (5 , 1 ) rectangle (0,0);
\fill [green!70] (4 , 2 ) rectangle (5,1);
}
  \foreach \x in {0,...,5}
      \draw (\x,0)--(\x,5) (0,\x)--(5,\x) ;
\end{tikzpicture}
      {\bf {VI}}
     \end{minipage}
     \caption{Illustration of the hybrid algorithm for a pixel grid. Each image represents one step of the algorithm, grey cells are yet to be updated. The star pixels are pixels in $\mathcal{D}$ with $A_{i,j}$ smallest. Green represents the update by GSLEX-I and the yellow pixels are the partial line smoothing updates. \label{fig:illusalg}}
  \end{figure}

In this new strategy, the only technical issue to address is that, at a pixel $(i, j)$ in set $\mathcal{D}$, the lagged
coefficient (here $A_{i,j}$) must be a previously updated pixel in this iteration otherwise we cannot avoid
multiple updates (as with Hybrid Smoother 1) within one smoothing iteration.
Our proposed solution is to view a group of adjacent pixels in set 
$\mathcal{D}$ whose smallest coefficient is $A_{i,j}$ (shown as starred pixels in Figure~\ref{fig:illusalg}) and sit on a line as a superpixel and to update together with their $A_{i,j}$ terms lagged. If the superpixel is comprised of a single pixel, we set its immediate neighbour pixel
(here $(i, j + 1)$) as a starred pixel so the group is of size 2. All other pixels in set $\mathcal{D}$ (without smallest coefficient $A_{i,j}$) and those not in $\mathcal{D}$ are treated as
normal pixels (non-starred) and are relaxed by the GSLEX-1 formula. Hence in each smoothing step,
starred and non-starred pixels are only updated once. 

 In Figure~\ref{fig:illusalg} we illustrate how this proposed algorithm would update the pixels, steps {\bf I--VI} represent one iteration of the smoother on the $5\times 5$ grid. The starred pixels represent those pixels which have $A_{i,j}$ the smallest. The algorithm proceeds as follows:
\begin{enumerate}
\item[{\bf I\ \ }]{We identify the pixels in $\mathcal{D}$ which have $A_{i,j}$ the smallest (indicated by a star).}
\item[{\bf II\ }]{Perform GSLEX-I on all non-starred pixels.}
\item[{\bf III}]{Collective partial line update on adjacent starred pixels.}
\item[{\bf IV}]{Perform GSLEX-I again on all non-starred pixels.}
\item[{\bf V\ }]{If a single starred pixel is found, update collectively with the immediate neighbour.}
\item[{\bf VI}]{Perform GSLEX-I again on all non-starred pixels.}
\end{enumerate}

\subsubsection{The adapted iterative schemes for other cases}

We previously focussed on the case for $A_{i,j}$ being lagged and now discuss other components of our
iterative scheme to cover the cases of $B_{i,j},C_{i,j}$ and $D_{i,j}$ being lagged. 

Crucially, to ensure that the scheme agrees with the LFA we must change the direction of update between the schemes for updating $A_{i,j}$, $B_{i,j}$, $C_{i,j}$ and $D_{i,j}$. For example, if we are lagging $B_{i,j}$ pixels we must update from the bottom-right corner to the top-left moving along rows right to left and from the bottom row to the top row. In Figure~\ref{fig:updatedir} we show the order in which the pixels should be updated for each lagging.

These sweeps in other directions are required to help those pixels in $\mathcal{D}$ that were treated as non-starred
pixels due to their smallest coefficients not being considered in the other sweeps. That is to say, each of 4
sweeps takes care of one type of alignment of the smallest coefficients (of course there are no other
directions to consider). Consequently, after all 4 sweeps, the compounded smoothing rate at each
pixel is small because we have ensured that one of the four multiplying factors is small while the
other three are no more than 1.

The broad algorithm ({\bf I--VI}) is the same in these cases as for the case of $A_{i,j}$ lagged; we identify the pixels which are of that case, perform GSLEX-I on all others and partial line updates on identified pixels.

Hybrid Smoother 2 performs 4 sweeps of the grid, each repeating the above {\bf I--V} and differing only in update order and assignment of starred pixels. In Figure~\ref{fig:updatedir} we display the order in which the pixels and superpixels should be updated for each lagging.

\begin{figure}[htb!]
  \begin{minipage}{0.15\linewidth}
    \centering
      \begin{tikzpicture}[darkstyle/.style={circle,draw,fill=black,minimum size=1,scale = 0.5},scale = 0.4]
      \foreach \x in {0,...,2}
   \foreach \y in {0,...,2}
{
\fill [gray!30] (0, 0 ) rectangle (5,5);
}
  \foreach \x in {0,...,5}
  {
      \draw (\x,0)--(\x,5) (0,\x)--(5,\x) ;
      }
        \foreach \x in {1,...,5}
        \foreach \y in {1,...,5}
  {
      \pgfmathtruncatemacro{\label}{(\x) +  ((5-\y))*5 }
            \node at ( \x-0.5,\y -0.5) {\label};
            }

\end{tikzpicture}

      {\bf {${\bm A_{i,j}}$ Lagged}}
  \end{minipage}
      \begin{minipage}{0.15\linewidth}
    \centering
  \begin{tikzpicture}[darkstyle/.style={circle,draw,fill=black,minimum size=1,scale = 0.5},scale = 0.4]
      \foreach \x in {0,...,2}
   \foreach \y in {0,...,2}
{
\fill [gray!30] (0, 0 ) rectangle (5,5);
}
  \foreach \x in {0,...,5}
      \draw (\x,0)--(\x,5) (0,\x)--(5,\x) ;
              \foreach \x in {1,...,5}
        \foreach \y in {1,...,5}
  {
      \pgfmathtruncatemacro{\label}{6-\x +  (\y-1)*5 }
            \node at ( \x-0.5,\y -0.5) {\label};
            }
\end{tikzpicture}
  {\bf {${\bm B_{i,j}}$ Lagged}}
     \end{minipage}
         \begin{minipage}{0.15\linewidth}
    \centering
    \begin{tikzpicture}[darkstyle/.style={circle,draw,fill=black,minimum size=1,scale = 0.5},scale = 0.4]
         \foreach \x in {0,...,2}
   \foreach \y in {0,...,2}
{
\fill [gray!30] (0, 0 ) rectangle (5,5);
}
  \foreach \x in {0,...,5}
      \draw (\x,0)--(\x,5) (0,\x)--(5,\x) ;
              \foreach \x in {1,...,5}
        \foreach \y in {1,...,5}
  {
       \pgfmathtruncatemacro{\label}{6-\y +  (\x-1)*5 }
            \node at ( \x-0.5,\y -0.5) {\label};
            }
\end{tikzpicture}
      {\bf {${\bm C_{i,j}}$ Lagged}}
     \end{minipage}
         \begin{minipage}{0.15\linewidth}
    \centering
    \begin{tikzpicture}[darkstyle/.style={circle,draw,fill=black,minimum size=1,scale = 0.5},scale = 0.4]
      \foreach \x in {0,...,2}
   \foreach \y in {0,...,2}
{
\fill [gray!30] (0, 0 ) rectangle (5,5);
}
  \foreach \x in {0,...,5}
      \draw (\x,0)--(\x,5) (0,\x)--(5,\x) ;
              \foreach \x in {1,...,5}
        \foreach \y in {1,...,5}
  {
      \pgfmathtruncatemacro{\label}{\y +  (6-\x-1)*5 }
            \node at ( \x-0.5,\y -0.5) {\label};
            }
\end{tikzpicture}
      {\bf {${\bm D_{i,j}}$ Lagged}}
     \end{minipage}
     \caption{Illustration of the hybrid algorithm for a pixel grid. The star pixels are pixels in $\mathcal{D}$ with $A_{i,j}$ smallest. Green represents the update by GSLEX-I and the yellow pixels are the partial line smoothing updates. \label{fig:updatedir}}
  \end{figure}

\subsubsection{Implementing Hybrid Smoother 2}

To ensure all laggings are considered, we sweep for $A_{i,j}$, $B_{i,j}$, $C_{i,j}$ and $D_{i,j}$ in this order, performing steps ({\bf I--VI}) on each sweep. These schemes are performed on all pixels in $\mathcal{D}$ and we see from Table~\ref{tab:comparerates3} that the maximum smoothing rate over $\mathcal{D}$ falls from 0.9997 to 0.4789. Therefore to reduce high-frequency errors by 90\%, with GSLINE-I this would have needed 7675 iterations but with the adapted iterative schemes we need only 4.

To ensure that all cases are considered, we design a hybrid smoother for which one outer iteration includes four sweeps of the image domain. In the first sweep we lag $A_{i,j}$, then in the second $B_{i,j}$ and so on. We note, for example, that in the sweep with $A_{i,j}$ lagged, then the pixels with coefficient $B_{i,j}$ smallest have a poor smoothing rate, however on the $B_{i,j}$ sweep the rate is good for these pixels and poor for those where we have  $A_{i,j}$ smallest. However, as the effects compound multiplicatively, after each outer iteration, the smoothing rate at pixels in $\mathcal{D}$ is good and for $\Omega\backslash\mathcal{D}$ is also good as these have had 4 GSLEX-I iterations.

We now consider the smoothing rates we can attain with this smoother. Firstly, for the Rada-Chen model \cite{SC28}, using (\ref{tab:compnonD}) we see that the maximum smoothing rate in each outer iteration of the smoother on $\Omega\backslash\mathcal{D}$ is approximately $0.7705^{4} = 0.3524$. By performing the adapted iterative schemes on $\mathcal{D}$
we have a maximum smoothing rate of 0.4789 (Table~\ref{tab:comparerates3}) in a single sweep. We know that the rate for GSLEX-I is poor for these pixels in $\mathcal{D}$ (close to 1) so the main reduction in error occurs when we perform the adapted scheme with the appropriate lagging. Therefore the maximum smoothing rate in one outer iteration of the smoother is approximately 0.4789, which is very good. 
 One consideration we must make is that the domain is covered $4$ times in each outer iteration, which could be computationally intensive for a large number of smoothing steps. Typically we find that for non-linear problems the number of overall sweeps of the grid is around 10-20 (see, for example, \cite{CB,SC31}) for the smoother, therefore we suggest 2 outer iterations (8 grid sweeps) which gives an impressive smoothing rate and is acceptable computationally.

{\bf Adaptive iterative schemes applied to the Spencer-Chen model} \cite{SC}.
We applied Hybrid Smoother $2$ to the Spencer-Chen model. In this case using just GSLINE-I
we have a maximum smoothing rate of 0.9990 but using the new smoother,
the maximum smoothing rate falls to 0.5032. Therefore, to reduce errors by 90\%
we need $4$ iterations rather than $2302$. This is a further indication that the
technique of using the partial line smoothers at the pixels with jumps in the
coefficients is a good way to reduce the maximum smoothing rate of the smoother
and the idea transfers to other models.

{\bf Improved smoothing rates for other images.} We now show how the maximum smoothing rate for Hybrid Smoother $2$ is smaller than GSLINE-I for several images with different levels of Gaussian noise. We compare to GSLINE-I as this is the recommended standard smoother for problems with jumping coefficients. We denote the corresponding maximum smoothing rates as $\mu_{GSLINE-I}$ and $\mu_{GSHYBRID}$ respectively. Results obtained previously are just for the clean image in Figure~\ref{fig:abcdcompare}(a). Here we compare the smoothing rates for noisy versions of this image and also of those in Figure~\ref{fig_6}.
\begin{table}[htb!]
\resizebox{\columnwidth}{!}{
\begin{minipage}{0.5\textwidth}
\begin{tabular}{l|c|c}
Image & $\mu_{GSLINE-I}$ & $\mu_{HYBRID}$  \\
\hline
\vspace{-0.4cm}&&\\
Figure~\ref{fig:abcdcompare}(a) + 1\% Noise 	& 0.9743 & 0.4891 \\
Figure~\ref{fig:abcdcompare}(a) + 5\% Noise 	& 0.9851 & 0.4815 \\
Problem 1 				& 0.9960 & 0.4532 \\
Problem 1 + 1\% Noise	& 0.9900 & 0.4749 \\
Problem 1 + 5\% Noise	& 0.9991 & 0.4789 \\
\end{tabular}
\end{minipage}\qquad 
\begin{minipage}{0.5\textwidth}
\begin{tabular}{l|c|c}
Image & $\mu_{GSLINE-I}$ & $\mu_{HYBRID}$  \\
\hline
\vspace{-0.4cm}&&\\
Problem 2 			& 0.9999 & 0.4736 \\
Problem 2 + 1\% Noise	& 0.9988 & 0.4886 \\
Problem 2 + 5\% Noise	& 0.9934 & 0.4518 \\
Problem 3 			& 0.9999 & 0.4829 \\
Problem 3 + 1\% Noise	& 0.9999 & 0.4863 \\
Problem 3 + 5\% Noise	& 0.9999 & 0.4841 \\
\end{tabular}
\end{minipage}
}
\caption{Comparison of the maximum smoothing rates for GSLINE-I and Hybrid Smoother $2$ for various images. \label{tab:comparerates33}}
\end{table}

{\bf Algorithm 3.} In future discussion, we refer to the Full Approximation Scheme using Hybrid Smoother 2 as Algorithm 3.

%
%
%


\section{Numerical Experiments}\label{sec:numerical}

\begin{figure}[htb!]
  \begin{minipage}{0.32\linewidth}
    \centering
      \includegraphics[width=\linewidth]{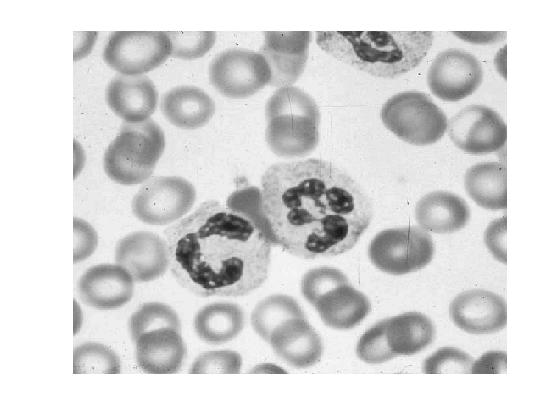}\\
       Problem 1
  \end{minipage}
  \begin{minipage}{0.32\linewidth}
    \centering
       \includegraphics[width=\linewidth]{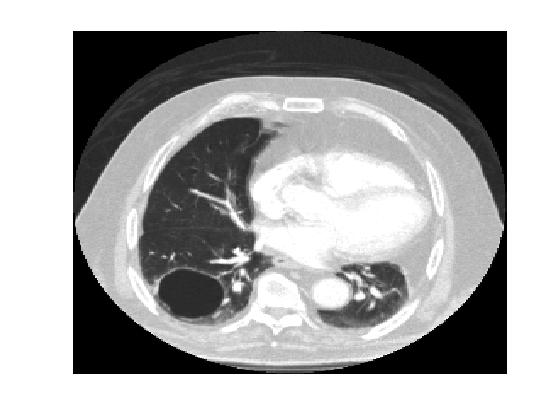}\\
       Problem 2
         \end{minipage}
    \begin{minipage}{0.32\linewidth}
    \centering
       \includegraphics[width=\linewidth]{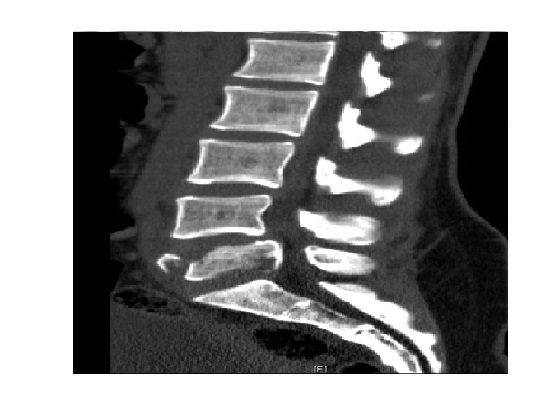}\\
       Problem 3
  \end{minipage}
\caption{The test images used in this section for the experiments.}
\label{fig_6}
\end{figure}

In this section we show two types of numerical experiments:
comparisons with the current best methods and analysis of the complexity of
Algorithms 2 and 3. Results have been obtained for many artificial and real images but we restrict to
the images shown in Figure~\ref{fig_6}.
We show real images as these are of most interest for the application
of selective segmentation.
The Rada-Chen and Spencer-Chen models we look at are non-convex and
we therefore need the initialisation to be close to the final solution.
Thankfully this can be achieved by setting the initial contour as the boundary of
the polygon formed from the user selected points in $\mathcal{S}$.
For examples of such user defined points, see Figure~\ref{fig:test11}.

{\bf Parameter Choices.} The values of $c_{1}$ and $c_{2}$,
being the average intensities inside and outside of the contour, are updated at the
end of each multigrid iteration - the initial values are set to the average inside and
outside the initial contour. We fix $\mu = 1/2$, $\lambda_{1} = \lambda_{2} = 10^{-4}$, $\nu = 1$ (for the Rada-Chen model) and $\theta = 1$ (for the Spencer-Chen model). In all experiments we use a V-cycle, i.e. fix $\gamma = 1$.

\begin{figure}[htb!]
\includegraphics[height=4cm,width=0.4\textwidth]{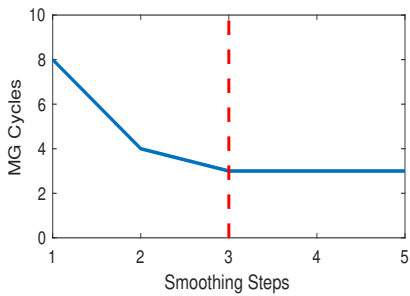}
\caption{The number of smoothing steps plotted against the number of multigrid cycles required to achieve convergence for Algorithm 3 on Problem 1. Guided by this, we choose 3 smoothing steps as the gain plateau's at this point. \label{fig:alg3sm}}
\end{figure}

{\bf Number of Smoothing Steps.}
To decide how many smoothing steps were required in Algorithms 1, 2 and 3, we performed experiments to see how the number of smoothing steps impacted the number of multigrid cycles for convergence. As the number of smoothing steps increases, the number of cycles decreases and plateaus. We fix the number of smoothing steps for each algorithm as the number required for the number of multigrid cycles to first plateau. In Figure~\ref{fig:alg3sm} we demonstrate how the number of multigrid cycles required for convergence changes with the number of smoothing steps and how we choose the optimal number of pre- and post-smoothing steps ($\nu_{1}$ and $\nu_{2}$). In all tests we use 100 iterations
of the exact solver (AOS) on the coarsest level.
Using this technique, we fix the smoothing steps for Algorithms 1, 2 and 3 as $\nu_{1} = \nu_{2} =$ 5, 3 and 3 respectively.

%
%

\subsection{Comparison of Algorithm $2$ and Algorithm $3$ with AOS}
\label{sec:alg23}

In this section we compare the speed of the proposed Algorithms 2 and 3 with AOS. We use the image from Problem 1 and scale this to different resolutions. The methods both use the standard stopping criteria
$\frac{||\phi^{(k+1)}-\phi^{(k)}||_{2}}{||\phi^{(k)}||_{2}}< \eta$,
where $\eta$ is a small tolerance parameter. In Table~\ref{tab:comparetimes} we see that Algorithm $3$ is faster to reach the stopping criteria (with $\eta = 10^{-4}$) than Algorithm $2$ and that both are faster than AOS for all but the smallest resolution image.
We see that as the image size grows larger, performance is significantly better.
One key aspect of Algorithms $2$ and $3$ is that we have the expected ratio for an $\mathcal{O}(N)$ method (in 2D) of 4 and hence an optimal complexity multigrid method. We also see that the multigrid method has a stable number of overall iterations, whereas with the AOS method, the iteration number grows as the image size grows.
Finally, we see that, although it converges faster overall, the cost per MG cycle is larger for Algorithm $3$ than $2$. This is due to a higher number of grid sweeps being required in the smoothing steps, however we believe that with improved and optimised coding of the smoother the performance of Algorithm $3$ can be increased to achieve far faster convergence than that of Algorithm $2$.
\begin{table}[htb!]
\resizebox{0.9\columnwidth}{!}{
\begin{tabular}{||c|c|c|c|c|c|c|c|c|c||}
\hline
\rule{0pt}{2.5ex}
 \multirow{2}{*}{Image size} &  \multirow{2}{*}{\begin{minipage}{2.2cm}
 \centering Number of Unknowns, $N$\end{minipage}}  &
 \multicolumn{2}{c|}{AOS} &  \multicolumn{3}{c|}{Algorithm 2} & \multicolumn{3}{c||}{Algorithm 3} 
 \\
\cline{3-10}
\rule{0pt}{2.5ex}
 & & Iter & CPU Time (s) & Iter & CPU Time (s) & CPU Ratio & Iter & CPU Time (s) & CPU Ratio \\
\hline
\vspace{-0.4cm}&&&&&&&&&\\
256 $\times$ 256    		&65536		& 32  	& 3.2	& 4	& 3.1		& - 		& 4& 8.8		& -	 \\
512 $\times$ 512  		&262144		& 39 	& 17.3 	& 5	& 11.6	& 3.7   	& 3&15.0		& 1.7\\
1024 $\times$ 1024 		&1048576		& 48 	& 123.5	& 5	& 44.0	& 3.8	& 3&	43.8		& 2.9\\
2048 $\times$ 2048 	&4194304	& 60 	& 759.2	& 5	& 174.2	& 4.0	& 3&	174.1		& 4.0\\
4096 $\times$ 4096 	&16777216	& 75 		& 8632.4	& 5	& 725.9	& 4.2 	& 3&	688.2	& 4.0\\
8192 $\times$ 8192 		&67108864	& *		& *		& 5	& 2952.2	& 4.1		& 3&	2766.9	& 4.0\\
\hline
\end{tabular}
}
\caption{For an image of size $N=m\times n$, we show a comparison of the number of iterations and the associated CPU times to achieve the same results for the Rada-Chen model for AOS and Algorithms $2$ and $3$. `*' indicates that the runtime exceeded 24 hours.}\label{tab:comparetimes}
\end{table}

\vspace{-0.2in}
\subsection{Comparison of Algorithms 1, 2 and 3}

We now look to see the practical gains from improving the smoother, i.e. the improved smoothing rate of Algorithm 3 should translate into a faster convergence rate \cite{Napov2011}.

\begin{defn}
In both Algorithms 2 and 3 we must identify the set $\mathcal{D}$, being pixels at which the coefficients vary significantly. To do this we compute the minimum multiplicative factor between the largest and smallest of the coefficients $A_{i,j},B_{i,j},C_{i,j},D_{i,j}$ (see \S\ref{sec:adaptiveschemes}). We will denote the minimum multiplicative factor by $\Sigma$.
\end{defn}

For completion, we will compare Algorithms 2 and 3 to Algorithm 1 for a range of $\Sigma$ values. The algorithms are all used to segment the image in Figure~\ref{fig:abcdcompare}(a), with fine grid $1024^{2}$ and coarse grid $32^{2}$ and $\eta=10^{-4}$ (all parameters are as earlier in \S\ref{sec:numerical}).

{\bf Level set energies.} In Table~\ref{tab:compareenergies} we give the energy of the level set at the end of each multigrid cycle for the Rada-Chen model for Algorithms 1, 2 and 3 for various $\Sigma$ values. The rows are ordered in descending order.
\begin{table}[htb!]
\begin{tabular}{||c|c|c|c|c|c|c|c||}
\hline
\rule{0pt}{2.5ex}
 \multirow{2}{*}{} &
 \multicolumn{7}{c||}{Iteration} \\
\cline{2-8}
\rule{0pt}{2.5ex}
& 1 & 2 & 3 & 4 & 5 & 6 & 7 \\
\hline
\vspace{-0.4cm}&&&&&&&\\
Algorithm 1     					& 2.4687& 1.9333& 1.9271& 1.9253 	& 1.9247 & 1.9241 	& 1.9236\\
Algorithm 2 ($\Sigma = 16$) 		& 2.4684& 1.9333& 1.9264& 1.9244 & 1.9238 &  - 		&  -\\
--------"-------- ($\Sigma = 8$)  	& 2.4683& 1.9321& 1.9251& 1.9242 	& 1.9235 &  - 		&  -\\
--------"-------- ($\Sigma = 4$) 	& 2.4683& 1.9302& 1.9242& 1.9237 & 1.9226  &  - 		&  - \\
--------"-------- ($\Sigma = 2$)  	& 2.4563 & 1.9269& 1.9214& 1.9207 & 1.9199  &  - 		&  - \\
Algorithm 3 ($\Sigma = 16$)		& 2.4300& 1.9185& 1.9180 &  - 	&  - 		&  - 		&  -  \\
--------"-------- ($\Sigma = 8$)  	& 2.4253& 1.9171& 1.9166 &  - 		&  - 		&  - 		&  -  \\
--------"-------- ($\Sigma = 4$)  	& 2.4184& 1.9167& 1.9164 &  - 		&  - 		&  - 		&  -  \\
--------"-------- ($\Sigma = 2$)  	& 2.4136& 1.9165& 1.9163 &  - 		&  - 		&  -  		&  -  \\
\hline
\end{tabular}
\caption{Level set energies ($\times 10^{5}$) after each multigrid iteration of Algorithms 1, 2 and 3 (for varying $\Sigma$) on the image in Figure~\ref{fig:abcdcompare}(a) + 10\% Gaussian noise. A dash indicates convergence before iteration number was reached.}\label{tab:compareenergies}
\end{table}

Firstly, we see that Algorithm 3 converges in $3$ cycles, where Algorithm 2 converges in $5$ and Algorithm 1 converges in $7$ cycles. Secondly, we notice that the energy is smallest for Algorithm $3$ and Algorithm $2$ gives a lower energy than Algorithm $1$ (for all $\Sigma$ values). Finally, we notice that as $\Sigma$ gets smaller (and the number of pixels in $\mathcal{D}$ increases), the energy of the level set at each cycle is smaller. This is all in agreement with the theoretical understanding of the smoothers, that they should give a small rate on the pixels in $\mathcal{D}$, and by increasing the size of $\mathcal{D}$ convergence improves.

{\bf Recommended Algorithm.} The CPU timings for Algorithm 3 are the best of the three algorithms (Table~\ref{tab:comparetimes}). The level set energies are also the lowest for Algorithm 3 (Table~\ref{tab:compareenergies}) at each iteration. It performs the best at tackling the PDEs which have many discontinuous coefficients and the experimental results are in agreement with the theory in \S\ref{sec:hyb1}. We therefore recommend Algorithm $3$ to achieve a fast solution to the Rada-Chen and Spencer-Chen selective segmentation models.
%
%

\begin{figure}[htb!]
  \begin{minipage}[b]{0.24\linewidth}
    \centering
\includegraphics[width=1\textwidth]{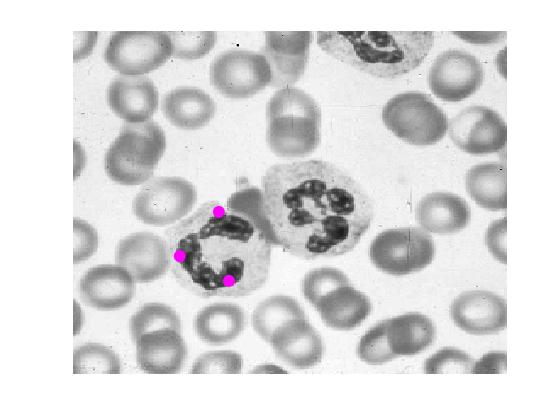}
  \end{minipage}
   \begin{minipage}[b]{0.24\linewidth}
    \centering
\includegraphics[width=\textwidth]{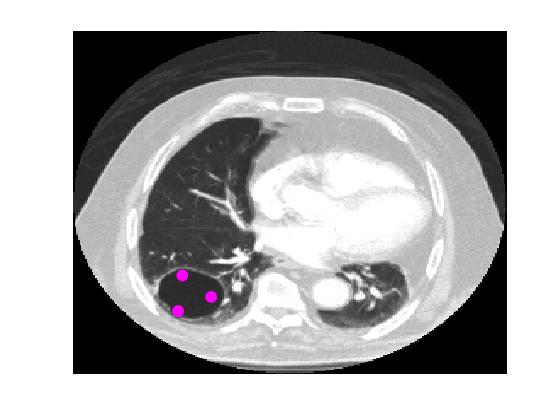}
  \end{minipage}
 \begin{minipage}[b]{0.24\linewidth}
    \centering
\includegraphics[width=\textwidth]{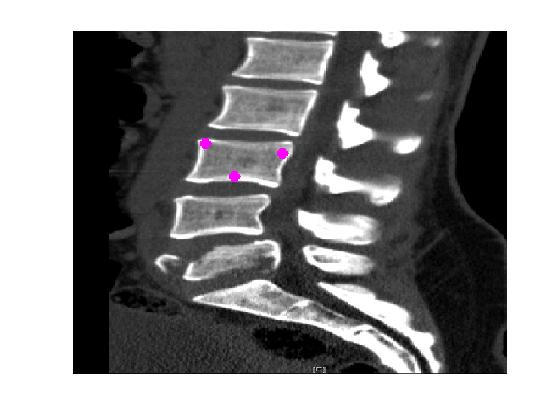}
  \end{minipage}
  \begin{minipage}[b]{0.24\linewidth}
    \centering
\includegraphics[width=\textwidth]{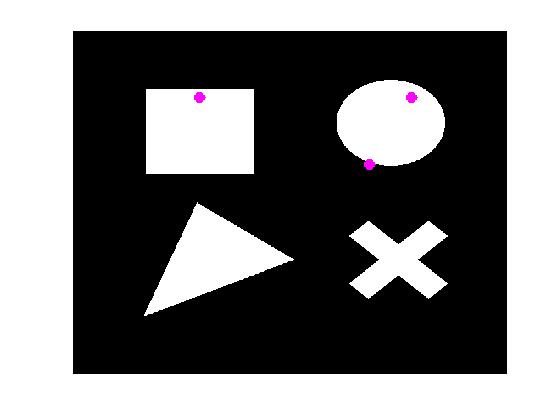}
\end{minipage}\qquad%
  \begin{minipage}[b]{0.24\linewidth}
    \centering
\includegraphics[width=1\textwidth]{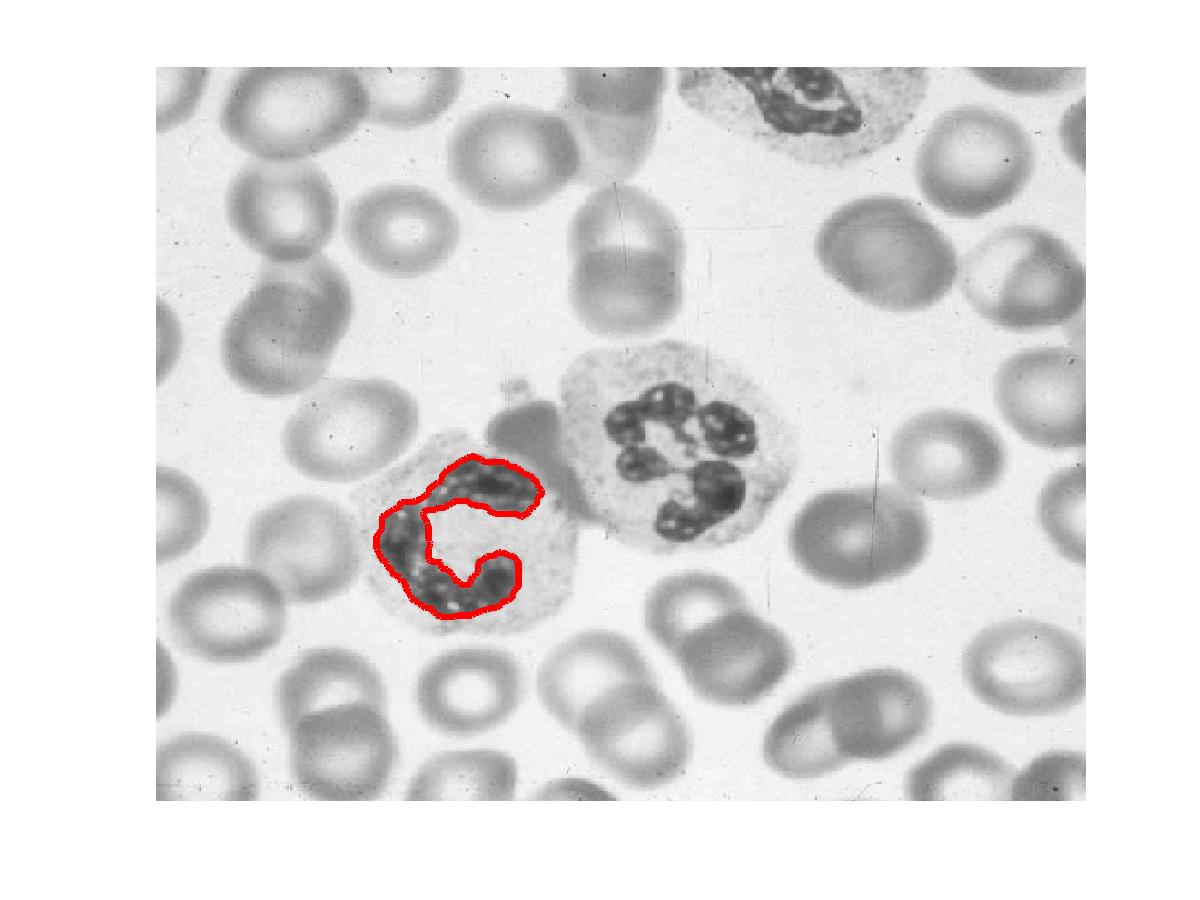}
  \end{minipage}
  \begin{minipage}[b]{0.24\linewidth}
    \centering
\includegraphics[width=\textwidth]{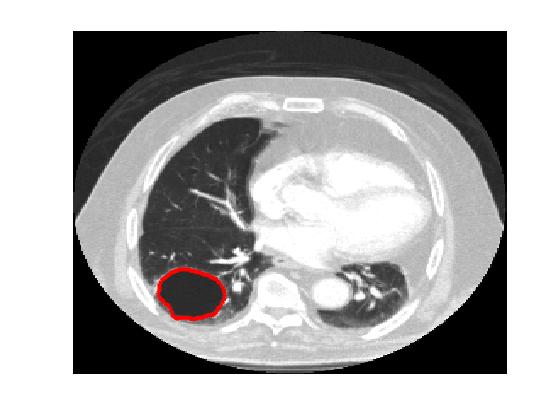}
  \end{minipage}
  \begin{minipage}[b]{0.24\linewidth}
    \centering
\includegraphics[width=\textwidth]{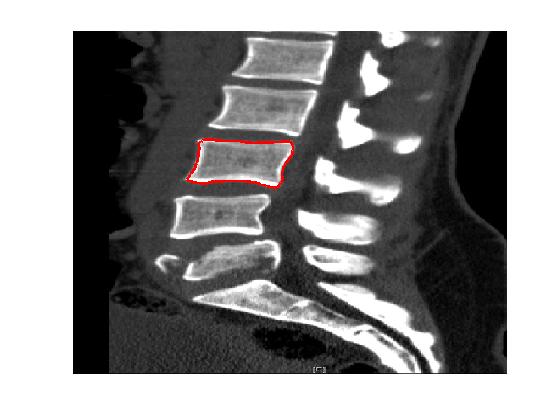}
  \end{minipage}
  \begin{minipage}[b]{0.24\linewidth}
    \centering
\includegraphics[width=\textwidth]{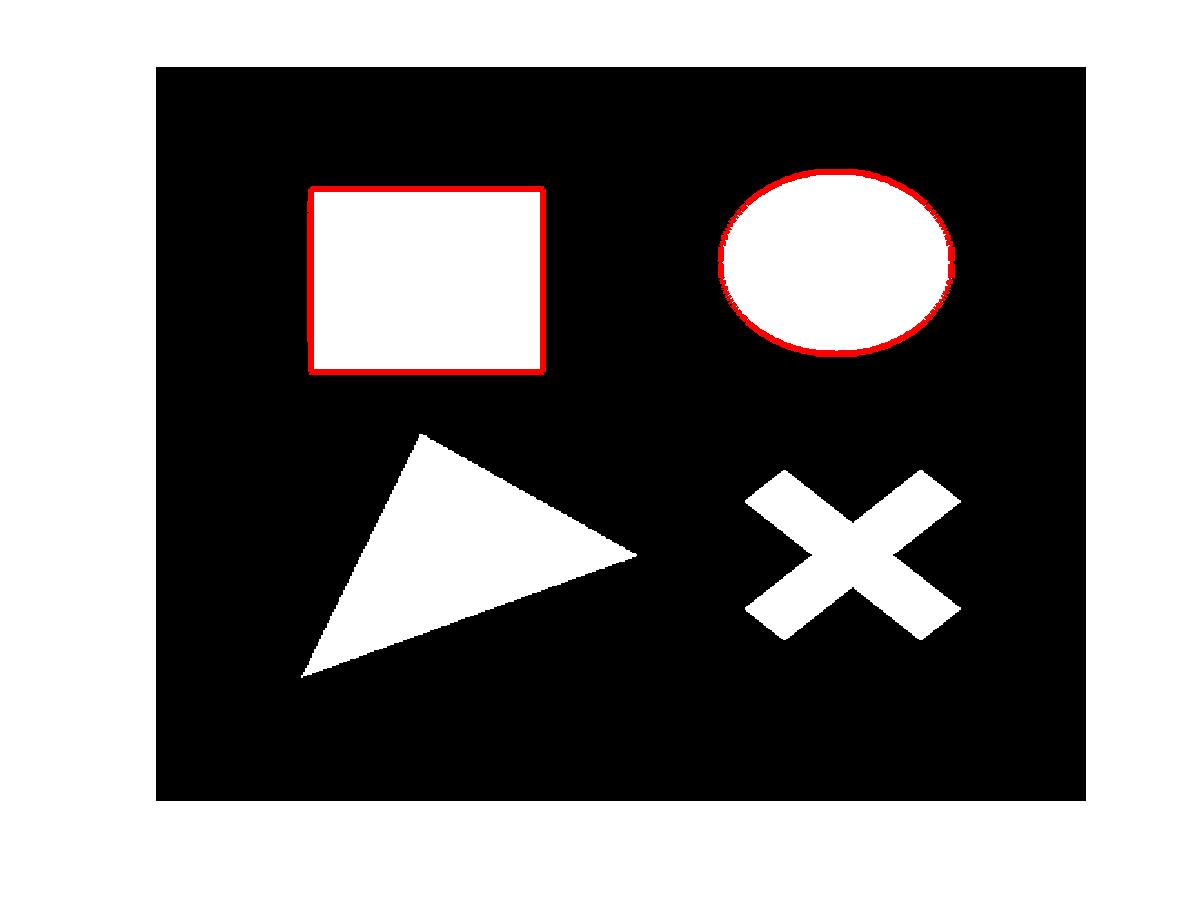}
\end{minipage}
\caption{Algorithm 3 results; user selections and segmentation results.}
\label{fig:test11}
\end{figure}

{\bf Algorithm 3 Results.} In Figure~\ref{fig:test11} we briefly show the results of
Algorithm $3$ applied to the test images for the Rada-Chen model shown in Figure~\ref{fig:abcdcompare}(a) and Figure~\ref{fig_6} with $\eta=10^{-4}$.
\vspace{-0.1in}

\subsection{Complexity of Algorithm 3}
We analyse Algorithm 3 to estimate the complexity of each multigrid cycle. We show analytically and experimentally that Algorithm 3 is $\mathcal{O}(N)$ as is expected for a multigrid method. We start with analysis of the complexity of the smoother, restriction operator, interpolation operator and coarse grid solver and then use the actual CPU times in Table~\ref{tab:comparetimes} to confirm the predicted complexity.

{\bf Analytical complexity.} Consider first only the fine grid with $N=nm$ pixels.
Hybrid smoother 2 uses GSLEX-I on $K$ pixels and partial line smoothers on $L$ segments, containing the remaining $N-K$ pixels. GSLEX-I requires 13K operations. The partial line smoothers require $\mathcal{O}(M_{i})$ operation, where $M_{i}$ is the size of the line segment for $i\in [0,L]$. Suppose the number of operations for each partial line smoothing is $\kappa M_{i}$. We can therefore bound the complexity of the smoothing as $13K + \kappa \sum_{i=0}^{L}M_{i}$.
We know that $K\le N$
and we perform $4$ grid sweeps for every $\nu_{1}$ pre-smoothing steps and $\nu_{2}$ post-smoothing steps. For simplicity, assume a square image (i.e. $n=m$) and so for smoothing on one level we have
\[
4(\nu_{1}+\nu_{2})\left(13K + \kappa \sum_{i=0}^{L}M_{i}\right) \le 4(\nu_{1}+\nu_{2})(13N + \kappa nL) \le 4(\nu_{1}+\nu_{2})(13+\kappa)N
\]
operations. With a $V$- cycle over $T$ grids, the number of operations is
\[
4(\nu_{1}+\nu_{2})(13+\kappa)N(1+\frac{1}{4}+\frac{1}{16}+\cdots +\frac{1}{2^{2(T-1)}}) < 4(13+\kappa)N \frac{\nu_{1}+\nu_{2}}{1-2^{-2}}=\frac{16(13+\kappa)(\nu_{1}+\nu_{2})}{3}N
\]
The restriction operator has complexity at most $15N$ on the finest grid and with $M$ grids there are $M-1$ restrictions, hence a complexity of less than $20N$. Interpolation has complexity at most $5N$ on the finest grid and hence all interpolation operators contribute at most $\frac{20}{3}N$ operations. Finally, with AOS as the coarse grid solver each iteration needs $448N\cdot 2^{-2(M-1)}$ operations, this is clearly bounded by $448N$. Therefore the overall maximum complexity of Algorithm $3$ is
\[
\frac{16(13+\kappa)(\nu_{1}+\nu_{2})}{3}N+ 20N + \frac{20}{3}N + 448\nu_{AOS}N \le \left[  \frac{16(13+\kappa)(\nu_{1}+\nu_{2})}{3}+448\nu_{AOS}\right] N,
\]
with $\nu_{AOS}$ the number of AOS iterations performed - as desired, the algorithm is $\mathcal{O}(N)$.

{\bf Experimental complexity.} In Table~\ref{tab:comparetimes} we show the ratio of the CPU times for Algorithm 3 on $\Omega^{h}$ when compared with the time on $\Omega^{2h}$. We see that the ratio is around 4 which linearly follows the increase in pixel number. Hence we see experimental confirmation of our analytical result that Algorithm 3 is an $\mathcal{O}(N)$ method.


\section{Conclusions}
Image segmentation models provide a set of challenging and non-linear
PDEs with non-smooth coefficients.
Direct application of multigrid solvers with standard smoothers such as the lexicographic and line Gauss-Seidel
smoothers
leads to poor or no convergence. This paper has investigated the reasons why smoothers
become ineffective due to non-smoothness of coefficients and proposed two hybrid smoothers
that are aware of jumps and add extra local smoothing using non-standard iterative schemes.
We find that both smoothers lead to convergent multigrid algorithms, however we recommend one smoother above the other as results are best experimentally and are shown to be good theoretically.
Experiments confirm that the proposed new algorithm,
outperforms the current fast methods. It also has optimal complexity and therefore
is suitable for solving selective segmentation models for large images.
Moreover, the ideas used in the design of the new smoother can be applied to
 other segmentation
models and potentially non-smooth PDEs from other applications.

\section*{Acknowledgements}
The first author wishes to thank the UK EPSRC and the Liverpool Heart and Chest Hospital for
supporting the work through an Industrial CASE award.
The second author is grateful to the
EPSRC for the grant EP/K036939/1.

\bibliography{biblio_V3}
\bibliographystyle{plain}
\end{document}